\documentclass[12pt]{article}
\usepackage{amssymb}

\usepackage{amsfonts}
\usepackage{amsmath}
\usepackage{amssymb,latexsym}
\usepackage[mathcal]{eucal}
\usepackage{amscd}

 \def\qed{\unskip\quad \hbox{\vrule\vbox to 6pt {\hrule width
          4pt\vfill\hrule}\vrule} }
          \newcommand{\bez}{\nopagebreak\hspace*{\fill}
          \nolinebreak$\qed$\vspace{5mm}\par}
          \newenvironment{proof}{\vspace{1ex}
               \vspace*{10mm}\vspace{-10mm}\par\noindent{\bf

               Proof.}\nopagebreak
                    \par}{\nopagebreak\linebreak[0]\hspace*{\fill}
 $\qed$\vspace{5mm}\pagebreak[0]\vspace{3ex}\par}

                         \newtheorem{Th}{Theorem}
                         \newtheorem{Prop}{Proposition}
                         \newtheorem{Lemma}{Lemma}
                         \newtheorem{Remark}{Remark}

                         \newtheorem{Cor}{Corollary}
                         \newtheorem{Problem}{Problem}

\newcommand{\al}{\alpha}

\newcommand{\cs}{{\cal S}}

\newcommand{\cc}{{\cal C}}

\newcommand{\ck}{{\cal K}}
\newcommand{\cb}{{\cal B}}
\newcommand{\ca}{{\cal A}}

\newcommand{\cj}{{\cal J}}

\newcommand{\cp}{{\cal P}}

\newcommand{\Q}{\mathbb{Q}}
\newcommand{\R}{{\mathbb{R}}}
\newcommand{\T}{{\mathbb{T}}}

\newcommand{\Z}{{\mathbb{Z}}}
\newcommand{\N}{{\mathbb{N}}}
\newcommand{\xbm}{(X,{\cal B},\mu)}

\newcommand{\ov}{\overline}
\newcommand{\beq}{\begin{equation}}
\newcommand{\eeq}{\end{equation}}
\newcommand{\vep}{\varepsilon}
\newcommand{\va}{\varphi}

\newcommand{\ot}{\otimes}
\newcommand{\la}{\lambda}
\newcommand{\wt}{\widetilde}

\newcommand{\supp}{{\rm{supp\,}}}
\newcommand{\Aut}{{\rm{Aut\,}}}
\newcommand{\Id}{{\rm{Id}}}

\newcommand{\Homeo}{{\rm{Homeo\,}}}

\begin{document}
\title{A topological lens for a measure-preserving system}
\author{Eli Glasner \and Mariusz Lema\'nczyk \and Benjamin Weiss}

\thanks{{\em 2000 Mathematical Subject Classification:
37A05, 37A35, 37B05, 37B40}}
\thanks{{\em Key words and phrases:} Topological lens, joinings, 
Devaney chaos, weak mixing, K-systems.}
\thanks{Research supported by the EU Program
Transfer of Knowledge ``Operator Theory Methods for Differential
Equations" TODEQ and the grant of Polish Ministry of Science and
High Education 2008-2012}

\date{}

\maketitle

\begin{abstract}
We introduce a functor which associates to every measure preserving
system $(X,\mathcal{B},\mu,T)$ a topological system 
$(C_2(\mu),\tilde{T})$ defined on the space of  2-fold couplings
of $\mu$, called the {\it topological lens} of $T$.
We show that often the topological lens
``magnifies" the basic measure dynamical properties of
$T$ in terms of the corresponding topological properties of $\tilde{T}$.
Some of our main results are as follows: (i) $T$ is weakly mixing iff
$\tilde{T}$ is topologically transitive (iff it is topologically weakly 
mixing). (ii) $T$ has zero entropy iff $\tilde{T}$ has zero
topological entropy, and $T$ has positive entropy iff $\tilde{T}$ has infinite topological entropy. (iii) For  $T$ a $K$-system, the topological lens is
a $P$-system (i.e. it is topologically transitive
and the set of periodic points is dense; such sytems are also
called chaotic in the sense of Devaney).
\end{abstract}

\newpage

\tableofcontents

\section*{Introduction}

Ergodic theory and topological dynamics are two branches of the
theory of dynamical systems. The first deals with groups acting on a
probability measure space in a measure-preserving way;
the second, with the action of groups on compact spaces
as groups of homeomorphisms.
Some of the terminology used in both branches is almost the same.
One speaks of transitivity, ergodicity, weak and strong mixing,
distality, rigidity, etc. both in ergodic theory and in topological
dynamics.
Even more surprising is the fact that major theorems in both
areas read almost the same. To mention one conspicuous
example, compare the statement of H. Furstenberg's theorem,
identifying topologically distal dynamical systems as inverse
limit of isometric extensions (\cite{Fu}),
with R. Zimmer's theorem, characterizing measure distal systems
(i.e. systems having  a separating sieve)
as systems admitting Furstenberg's towers
of (measure) isometric extensions (\cite{Z1},\cite{Z2}).

In the present work we restrict our attention to the classical
case where the acting group is the group of integers $\Z$.
We denote by $T$ the measure preserving transformation
(or the homeomorphism) which corresponds to $1 \in \Z$
and write $(X,\mathcal{B},\mu,T)$ (or $(X,T)$) for the general
measure preserving (or topological compact metric) system
we study.

Given a topological system $(X,T)$ we denote by $M(X)$
the compact convex set of probability measures on $X$
and by $M_T(X)$ the compact convex subset
(in fact the simplex)
of $T$-invariant measures. By a classical theorem
of Krylov and Bogoliouboff this is always a non-vacuous set.
When $M_T(X)$ consists of a single measure, say $M_T(X)=\{\mu\}$,
we say that $(X,T)$ is uniquely ergodic. In the general case,
each element $\mu$ of $M_T(X)$ defines a measure dynamical system
$(X,\mathcal{B}_X,\mu,T)$ ($\mathcal{B}_X$ denotes
the Borel $\sigma$-field of the topological space $X$).
We say that $\mu \in M(X)$ is {\it full} when
$\supp \mu =X$.

Part of the mystery of the elusive connection between
the two theories is removed by the following theorem
whose proof is straightforward;
one only has to note that for a full measure $\mu$,
$\mu(U) > 0$ whenever $U$ is a non-empty open set.

\begin{Th}
Let $(X,T)$ be a topological dynamical system,
$\mu\in M_T(X)$ a full measure;
then $(X,T)$ is topologically transitive, weakly mixing,
topologically mixing if the measure-preserving system
$(X,\mathcal{B}_X,\mu,T)$ is ergodic, weakly mixing,
mixing, respectively.
\end{Th}

A substantial part of modern ergodic theory deals with the
converse situation. One starts with a given measure ergodic system
$(X, \mathcal{B},\mu,T)$ and then looks for a {\it topological model};
i.e. a topological system $(Y,S)$ and a measure $\nu \in M_S(Y)$
such that the systems $(X,\mathcal{B},\mu,T)$ and $(Y, \mathcal{B}_Y,
\nu,S)$ are isomorphic, and such that the topological system $(Y,S)$
has some special properties like being minimal or uniquely ergodic.
The prototype for this kind of statements is the famous
Jewett- Krieger theorem which gaurantees the existence of a uniquely
ergodic model for any ergodic system.

In the present work we offer a novel perspective on the investigation
of the connection between measure and topological systems.
We introduce a functor which associates to every measure preserving
system $(X,\mathcal{B},\mu,T)$ a topological dynamical system
$(C_2(\mu),\tilde{T})$ defined on the space of  2-fold couplings
of $\mu$, called the {\it topological lens} of $T$.
More specifically, a {\em coupling} of $\mu$ is a probability measure $\xi$ on
$X \times X$ with both marginals equal to $\mu$.
We equip $C_2(\mu)$ with the (compact metrizable) weak$^*$ topology
and define $\tilde{T}: C_2(\mu) \to C_2(\mu)$ by the formula
$$
\tilde{T}(\xi)(A\times B) = \xi(T^{-1}A\times T^{-1}B),\ A,B \in
\mathcal{B},
$$
that is, $\tilde{T}(\xi)= (T\times T)_*(\xi)$.
As we will show, the topological lens
usually ``magnifies" the basic measure dynamical properties of
$T$ in terms of the corresponding topological properties of $\tilde{T}$.

Briefly our main results are as follows: (i) $T$ is weakly mixing iff
$\tilde{T}$ is topologically transitive (iff it is topologically weakly
mixing). (ii) $T$ has zero entropy iff $\tilde{T}$ has zero
topological entropy, and $T$ has positive entropy iff $\tilde{T}$ has infinite topological entropy.  (iii) $T$ is rigid iff $\tilde{T}$ is
pointwise recurrent (iff $\tilde{T}$ is uniformly rigid).
(iv) Distality of $\tilde{T}$ implies that $T$ has discrete spectrum.
(v) For the Bernoulli system $T$ of infinite entropy, the topological
lens $(C_2(\mu),\tilde{T})$ is a universal system both topologically
and measure theoretically in the sense that
every metric compact topologically transitive,
as well as every measure preserving system, appears as a subsystem of
$C_2(\mu)$. (vi) For  $T$ a $K$-system, the topological lens is
a $P$-system (i.e. it is topologically transitive
and the set of periodic points is dense; such systems are also
called chaotic in the sense of Devaney).
(vii) Finally, for many zero entropy measure systems (including
the generic automorphism $T \in \Aut(\mu)$) the set of
periodic points of the topological lens $\tilde{T}$ is closed
and nowhere-dense.

\section{The space of couplings}
Assume that $\xbm$ is a standard probability  Borel space. By
$C_2(X,\mu)$ (or $C_2(\mu)$) we denote the space of {\em
2-couplings} of $\xbm$, that is the space of probability measures
on $(X\times X,\cb\ot\cb)$ with projections $\mu$ on both
coordinates. The formula \beq\label{cor1-1} \langle
J_\rho(1_A),1_B\rangle_{L^2(\mu)}=\langle
1_A,1_B\rangle_{L^2(\rho)},\;\;A,B\in\cb\eeq
establishes a
one-to-one correspondence $\rho\mapsto J_\rho$ between $C_2(\mu)$
and the space $\cj(\mu)$ of doubly stochastic operators (Markov
operators) on $L^2\xbm$ ($J:L^2\xbm\to L^2\xbm$ is called {\em
doubly stochastic} if it is positive and $J(1)=J^\ast(1)=1$; note
that necessarily $\|J\| = 1$). With respect to the weak operator
topology $\cj(\mu)$ forms a compact semitopological
(metrizable) semigroup, where multiplication is
defined by composition: $J_{\rho_1\circ\rho_2}=J_{\rho_1}\circ
J_{\rho_2}$. Recall that a metric compatible with the weak
topology on $\cj(\mu)$ is given by 
\beq\label{compatible}
d(J,J')=\sum_{i,j=0}^\infty \frac1{2^{i+j}}|\langle
Jf_i,f_j\rangle -\langle J'f_i,f_j\rangle|, 
\eeq 
where $\{f_i\}_{i\geq0}$ is a
dense subset of $L^2\xbm$ and
$\|f_i\|=1$ for $i\geq0$. Using once more the correspondence
$\rho\leftrightarrow J_\rho$ we observe that in the weak topology
on $C_2(\mu)$ we have
$$
\rho_n\to\rho\;\;\mbox{iff}\;\;\rho_n(A\times B)\to\rho(A\times
B)\;\;\mbox{for all}\;A,B\in\cb.
$$
It also follows that a basis of
open sets in $C_2(\mu)$ is given by the family of sets of the form
\beq\label{postac} U(\alpha,\vep,P)=\{\rho\in C_2(\mu):\:
|\rho(A_i\times A_j)-p_{ij}|<\vep\},\eeq where
$\alpha=\{A_1,\ldots,A_k\}$ is a Borel partition of $X$, $\vep>0$
and
$P=[p_{ij}]_{i,j=1}^k$ is a non-negative matrix such that the
sum of elements in  the $i$-row is equal to $\mu(A_i)$ and the sum
of elements  in
the $j$-column is equal to $\mu(A_j)$, $1\leq
i,j\leq k$.

Denote by $\Aut(\mu)=\Aut\xbm$ the group of automorphisms of $\xbm$.
Notice that $\Aut\xbm$ naturally embeds into $C_2(\mu)$ where the 
embedding is given by
$S\mapsto\mu_S$.
Here $\mu_S$ stands for the graph measure given
by $S$, i.e. $\mu_S(A\times B)=\mu(S^{-1}A\cap B)$ for each
$A,B\in\cb$.
The Markov operator corresponding to $\mu_S$ is equal
to $U_S$ (the Koopman operator associated to $S$), where
$U_S(f)=f\circ S$ for each $f\in L^2\xbm$. Once $\mu$ is
understood we will also use
the notation $\Delta_S$ for $\mu_S$. The
embedding $S\mapsto \Delta_S$  is also topological, as on the
group of automorphisms considered as a (closed) subset of
$\mathcal{U}(L^2\xbm)$,
the group of unitary operators on $L^2\xbm$,
the weak and the strong operator topologies coincide.
It follows from a theorem of Kuratowski that the set of graph couplings,
that is the image of the embedding $\Aut(\mu)\subset C_2(\mu)$, is a
$G_\delta$ subset of $C_2(\mu)$. It is not hard to check that it is also
dense in $C_2(\mu)$ (see e.g. subsection \ref{Sec-WM} below).
We denote this dense $G_\delta$ subset of graph couplings
by $C_{gr}(\mu)$, but often we will regard $\Aut(\mu)$, via this
embedding, as a subset of $C_2(\mu)$.

Suppose $\ca\subset\cb$ is a $T$-invariant sub-$\sigma$-algebra. 
We can then
consider the factor system $(X/\ca,\ca,\mu)$ (instead of $\mu$ we
sometimes write  $\mu_{\ca}$ if we want to emphasize that we
consider the quotient system), here $X/\ca$ stands for classes of
points of $X$ that cannot be separated by sets of $\ca$. If
$\la\in C_2(\mu_{\ca})$ then we can lift it to an element of
$C_2(\mu)$, denoted by $\widehat{\la}$ and called the {\em
relatively independent extension of} $\la$, by setting
$$
\int_{X\times X}f\otimes g\,d
\widehat{\la}=\int_{X/\ca\times X/\ca}E(f|\ca)(\ov{x})
E(g|\ca)(\ov{y})\,d\la(\ov{x},\ov{y})
$$
for $f,g\in L^2\xbm$.
Denoting by $i_{L^2(\ca)}$ the natural embedding of $L^2(\ca)$ into
$L^2(\cb)$, by a simple calculation we have the following.

\begin{Lemma}\label{relind1}
$J_{\widehat{\la}}=i_{L^2(\ca)}\circ
J_{\la}\circ E(\cdot|\ca)$.\bez
\end{Lemma}

It follows immediately that the following holds.

\begin{Lemma}\label{relind2}
The map $C_2(\mu_{\ca})\ni\la\mapsto\widehat{\la}\in C_2(\mu)$ is
a continuous monomorphism of semitopological semigroups; i.e. it
is continuous, 1-1 and
$(\la\circ\rho)\,\widehat{}\,=\widehat{\la}\circ\widehat{\rho}$.
In particular this map is an embedding of the topological
system $C_2(\mu_A)$ into $C_2(\mu)$.
\bez
\end{Lemma}


Let $T$ be an element of $\Aut(\mu)$. By $J_2(T)\subset C_2(\mu)$ we
denote the set of {\em 2-self-joinings} of $T$. Recall that a
2-self-joining of $T$ is just
an element of $C_2(\mu)$ which is $T\times T$-invariant.
Equivalently, $\rho\in C_2(\mu)$ is a
$2$-self-joining of $T$ if and only if $J_\rho\circ U_T=U_T\circ
J_\rho$. From this it easily follows that $J_2(T)$ is a closed
subsemigroup of $C_2(\mu)$. When $T$ is ergodic $J^e_2(T)$, 
the set of ergodic $2$-self-joinings, is nonempty and it coincides with the 
collection of extremal points of the simplex $J_2(T)$. 
(Warning: $C_2(\mu)$ is not a simplex.) 
When $\ca\subset\cb$ is a $T$-invariant
sub-$\sigma$-algebra then the quotient action of $T$ on
$(X/\ca,\ca,\mu)$ is called a {\em factor} of $T$; we will often
write $T|_{\ca}$ to denote this action. It is easy to see that if
$\la\in J_2(T|_{\ca})$ then $\widehat{\la}\in J_2(T)$ and
therefore Lemma~\ref{relind2} is also true in the context of
self-joinings.

Extending the notion of $2$-self-joinings we define the set $J_n(T)$
($n\geq1$ or even $n=\infty$) of $n$-{\em self-joinings} of $T$.
These are the $T^{\times n}$-invariant probability measures
on $(X^n,\cb^{\ot n})$
all of whose one dimensional marginals are equal to $\mu$.

For a more about joinings we refer to \cite{Gl} and the list of references thereof.

\section{
The topological lens of an automorphism}\label{Sec-lens}
Given $T\in \Aut\,\xbm$ consider the
$\Z$-action $\wt{T}$
on $\cj(\mu)$ defined by conjugation:
\beq\label{wzor1} \wt{T}(J)=U_T^{-1}\circ J\circ U_T,\eeq where
$U_T$ is the Koopman operator associated to $T$. We will call this
action
the {\em topological lens of} $T$. Notice that
$\wt{T}$ is a homeomorphism of $\cj(\mu)$. Since
$$
\langle\wt{T}(J)1_A,1_B\rangle_{L^2(\mu)}=
\langle J\circ U_T1_A,U_T1_B\rangle_{L^2(\mu)}= \langle
J1_{T^{-1}A},1_{T^{-1}B}\rangle_{L^2(\mu)},
$$
the corresponding action on $C_2(\mu)$, which we also denote by
$\wt{T}$, is given by (see~(\ref{cor1-1}))\beq\label{wzor2}
\wt{T}(\rho)(A\times B)=\rho(T^{-1}A\times T^{-1}B)\;\;\mbox{for
each}\;A,B\in\cb.\eeq It is easy to see that
\beq\label{w1}\wt{T^n}=\wt{T}^n\eeq for each $n\in\Z$. 
Moreover if
$T_1$ (acting on $(X_1,\cb_1,\mu_1)$) is a factor of $T$, then
\beq\label{w2}\mbox{$\wt{T_1}$ is a topological factor of
$\wt{T}$}.\eeq 
In fact, if $\theta:\xbm\to(X_1,\cb_1,\mu_1)$
satisfies $\theta\circ T=T_1\circ\theta$, then the map
$\wt{\theta}(J)=V_\theta^\ast\circ J\circ V_\theta$
($V_\theta:L^2(X_1,\mu_1)\to L^2(X,\mu)$, $V_\theta(f_1)=f_1\circ
\theta)$ is
the corresponding continuous homomorphism of topological
dynamical systems.
Equivalently, if $\rho\in C_2(\mu)$ then
$$
\wt{\theta}(\rho)(A_1\times
B_1)=\rho(\theta^{-1}A_1\times\theta^{-1}B_1)$$ for each
$A_1,B_1\in\cb_1$. In order to show that $\wt{\theta}$ is onto we
use  the relative independent extension construction (see
Lemma~\ref{relind2} above).

Of course for every $T\in \Aut \xbm$ the
dense $G_\delta$ subset $C_{gr}(\mu)\subset C_2(\mu)$ is $\wt{T}$-invariant.
Since the action of $\wt{T}$ on $C_{gr}(\mu)$ is
isomorphic to conjugation by $T$ on $\Aut(\mu)$,
and as the group $\Aut\,(\mu)$ is algebraically
simple \cite{Fa}, it follows that the homomorphism $T \mapsto \wt{T}$
from $\Aut\,(\mu)$ into the group $\Homeo(C_2(\mu))$ is
an isomorphism. A more difficult question is whether
there are $S, T \in \Aut\,(\mu)$ which are not conjugate
(that is, the measure preserving systems  $(X,\mathcal{B},\mu,S)$
and $(X,\mathcal{B},\mu,T)$ are not isomorphic) while their topological
lenses $(C_2(\mu),\wt{S})$ and $(C_2(\mu),\wt{T})$ are isomorphic as
topological systems. It is not hard to see that this can not happen
for ergodic rotations (see Proposition \ref{DS} below). However,
it is very likely that the answer is yes.
In fact, it seems that the argument in the proof of Theorem
\ref{egmlbw} below might be refined to show that any two measure theoretical
Bernoulli systems have isomorphic lenses.

Notice that every (nontrivial) topological lens has many fixed points (in
particular,  $\wt{T}$ has many minimal subsets); indeed
the fixed point set of $\wt{T}$ coincides with the set of joinings,
\beq
\label{w3} Fix(\wt{T})=J_2(T).
\eeq
Notice that $J_2(T)$ is a
closed $\wt{T}$-invariant subset of $C_2(\mu)$ whose interior is
empty (to see this, given $\rho\in J_2(T)$ and $\vep>0$ take
first $\eta\in C_2(\mu)$ which is not a 2-self-joining and then
consider $\rho_{\vep}=(1-\vep)\rho+\vep\eta$ which cannot be a
2-self-joining for $\vep>0$). Similarly, periodic points for the
lens correspond to 2-self-joinings of powers of $T$. 
Thus the set of periodic points for
$\wt{T}$, which is just $\bigcup_{n\geq1}J_2(T^n)$,
is $\wt{T}$-invariant and meager.

The aim of this paper is to
investigate how topological properties of the lens
$(C_2(\mu),\wt{T})$ reflect ergodic properties of $T$. 
(Unless we say explicitly otherwise we usually assume that $T$ is ergodic.)

\vspace{0.3cm}

Before we begin our study of topological
lenses we briefly discuss  a
simpler topological system on $\cj(\mu)$ given by the one sided
composition $t_T$:
$$
t_T(J):=J\circ U_T\;\;\mbox{for}\;J\in \cj(\mu);$$ the
corresponding action on $C_2(\mu)$ which is given by the formula
$$
t_T(\rho)(A\times B)=\rho(T^{-1}A\times B)\;\;\mbox{for
each}\;A,B\in\cb
$$
we will also denote by $t_T$. We will now argue
that this system is particularly simple  from the dynamical point
of view. Indeed, because of
one sided continuity of the composition of
Markov operators (both left and right),
the enveloping semigroup of $t_T$ consists solely of
continuous maps $t_J$, where $J$ belongs to the weak closure of
the group generated by $U_T$. It follows that the topological
system  $(t_T,C_2(\mu))$ is weakly almost periodic (WAP); see
\cite{EN} and \cite{Gl}.
For such systems the closure of each orbit contains exactly
one minimal subsystem and this unique minimal system
is a compact monothetic topological group. 
In fact, we can describe these minimal (sub)systems quite precisely.

Denote by $\ck\subset\cb$ the sub-$\sigma$-algebra corresponding
to the Kronecker factor of $T$. Put
$$
\widehat{J^e_2(\ck)}=\{\widehat{\la}:\:\la\in
J^e_2(T|_{\ck},\mu_{\ck})\}
$$
and $\widehat{\cj^e(\ck)}$ for the corresponding set of Markov operators.

\begin{Prop}\label{charj}
There exists a sequence $(n_i)$ of density~1 such that for each
$\rho \in C_2(\mu)$ the limit points along
subsequences  $(n_{i_k})$ of
$t_{T}^{n_{i_k}}(\rho)$ are of the form
$\rho \circ\widehat{\la}$ for some $\la\in
J^e_2(T|_{\ck},\mu_{\ck})$.
\end{Prop}

\begin{proof}
The result holds trivially if the system $T$ has discrete
spectrum; indeed, recall that in this case the set of Markov
operators corresponding to  $J^e_2(T)$ is closed and equal to the
closure of powers of $U_T$. Take for $(n_i)$ the sequence of all
natural numbers.


Next assume that $T$ has partly continuous
spectrum. There exists a sequence $(n_i)$ of density~1 such that
if we denote \beq\label{decomp1} L^2\xbm=L^2(\ck)\oplus F\eeq then
for each $f,g\in F$
$$
\langle U_T^{n_i}f,g\rangle\to 0,\;\mbox{when}\;i\to\infty.
$$
Since $T|_{\ck}$ has discrete spectrum, we can choose a
subsequence $(n_{i_k})$ of $(n_i)$ such that
$U_{T|_{\ck}}^{n_{i_k}}\to V$, with $V\in \cj(T|_{\ck})$ and
$V=J_{\la}$ for some $\la\in J_2^e(T|_{\ck})$. All we need to show
is that
$$
U_{T}^{n_{i_k}}=t_{T}^{n_{i_k}}(Id)\to J_{\hat{\lambda}}.
$$
With no loss of generality we can assume $U_T^{n_{i_{k}}}\to J$.
Note that $J$ preserves the decomposition~(\ref{decomp1});
a weak limit of powers of $U_T$ preserves a weakly closed
$U_T$-invariant subspace.
For $f\in L^2\xbm$ and $g\in F$, by decomposing $f=f_1+f_2$, where
$f_1\in L^2(\ck)$ and $f_2\in F$,  we have
$$
\int_{X} f_2\circ
T^{n_{i_k}}\cdot g\,d\mu\to0
$$
by the property of $(n_i)$, hence
$$
\int_{X} f\circ T^{n_{i_k}}\cdot g\,d\mu\to\int_X J(f)\cdot
g\,d\mu=\int_XJ(f_1)\cdot g\,d\mu.
$$
Thus $Im(J)\subset L^2(\ck)$
or more precisely $J(L^2(\ck))\subset L^2(\ck)$ and
$J((L^2(\ck)^\perp)=\{0\}$. Hence, $J=i_{L^2(\ck)}\circ
J_\lambda\circ E(\cdot|\ck)=J_{\hat{\lambda}}$ by
Lemma~\ref{relind1}.
\end{proof}

Notice that the action of $t_{T|_{\ck}}$ on $J^e_2(T|_{\ck})$ just
goes back to the action of $T|_{\ck}$ on $(X/\ck,\ck,\mu)$ if we
recall that $T|_{\ck}$ is isomorphic to a minimal rotation on
a compact monothetic group $X/\ck$ (endowed with its normalized Haar
measure $\mu_{\ck}$). Therefore $t_T$ acting on
$\widehat{J^e_2(\ck)}$ is also minimal, and for each $J\in \cj(T)$
the only minimal subset contained in the orbit closure of $J$ is equal to
$J\circ\widehat{\cj^e(\ck)}$.

\begin{Remark}\label{attractor}\em
Let $S$ be a homeomorphism of a compact
metric space $M$. Let $A\subset M$ be closed. We say that $A$ is a
{\em quasi-attractor} if there exists a sequence $(n_i)$ of
density one such that for each $x\in M$, every limit point
of the sequence $S^{n_i}x$ lies in $A$. We claim
that in such a case: each invariant measure is concentrated on $A$
and each minimal subset is contained in $A$. In fact, suppose that
$\mu$ is ergodic for $S$ and let $x\in M$ be a generic point 
for $\mu$. Suppose that $\mu(M\setminus A)>0$. 
Choose a compact subset $C\subset M\setminus A$ for which $\mu(C)>0$ 
and then a continuous function $f$ with $0\leq f \leq1$ such that 
$f=1$ on $C$ and $f=0$ on $A$. Since $f=0$ on $A$ and $(n_i)$ 
has density~1,
$$
\frac1n\sum_{k=0}^{n-1}f(S^kx)\to 0.
$$ 
Hence $\int_Mf\,d\mu=0$, contradicting the 
fact that $\int_M f\, d\mu \ge \int_C f \, d\mu = \mu(C) > 0$.
\end{Remark}

We now observe that $J\circ\widehat{\cj^e(\ck)}$ is a
quasi-attractor in the closure of the orbit of $J\in \cj(\mu)$.
Finally note that if $T$ is weakly mixing, then all the sets
$J\circ\widehat{\cj^e(\ck)}$
collapse into one point
($\widehat{J^e_2(\ck)}=\{\mu\ot\mu\}$), and we obtain just one
attracting fixed point for the whole action $t_T$ on $C_2(\mu)$.

\vspace{0.3cm}

In contrast, the conjugation by $T \in \Aut(\mu)$ on the $\cj(\mu)$,
that is, the action of $\wt{T}$ on $C_2(\mu)$, is usually
not WAP. This will be amply demonstrated in the rest of this work,
but here is a simple first example.

\vspace{0.3cm}

{\bf Example.}\
Let $G=SL(2,\R)$ and $\Gamma\subset G$ a discrete cocompact
subgroup. The {\em geodesic} and {\em horocycle} flows on
the compact homogeneous space $X=G/\Gamma$
are the restrictions of the left multiplication
$G$-action  to the subgroups $\{A_t =
\left(
\begin{smallmatrix}
e^t & 0  \\  0 & e^{-t}
\end{smallmatrix}\right)
:t \in \R\}$
and $\{B_s = \left(
\begin{smallmatrix}
1 & s \\  0 &1
\end{smallmatrix}\right):s \in \R\}$,
respectively.
These flows preserve the normalized Haar measure $\lambda$ on $X$.
Now the commutation relations
$$
A_t B_s A_t^{-1}=
\left(\begin{smallmatrix}
1 & se^{2t} \\  0 &1
\end{smallmatrix}\right),
$$
applied to the corresponding unitary operators on $L_2(X,\la)$,
show that in the weak operator topology
\allowbreak $\lim_{t \to -\infty} \allowbreak A_t B_s A_t^{-1}= \Id$,
while, by mixing of the horocycle
flow, $\lim_{t \to \infty}\allowbreak A_t B_s A_t^{-1}= P$, where $P$ is the projection
onto the subspace of constant functions.
Thus in $\cj(\la)$ the orbit closure of, say $B_1$, under conjugation
by $A_1$ is isomorphic to the two-point compactification of $\Z$.
In particular, there are in this orbit closure two distinct minimal sets
$\{\Id\}$ and $\{P\}$, and therefore it is not WAP.

\section{Recurrence, distality and weak mixing}

In this section we impose certain topological conditions
on the lens like recurrence, distality and transitivity, and examine their
implication for the corresponding measure preserving transformation.

\subsection{Pointwise recurrence and rigidity}

Let us recall
some basic definitions. Let $R$ be a
homeomorphism of a compact space $Z$. Then $R$ is called {\em
pointwise recurrent} if for every $z\in Z$ and $\vep>0$ the orbit of
$z$ returns to the $\vep$-ball centered at $z$ infinitely often.
If there is an increasing sequence $(n_i)$ of integers such that
$R^{n_i}\to Id$ uniformly (pointwise) then $R$ is called {\em
uniformly rigid} ({\em rigid}).

An automorphism $T:\xbm\to\xbm$ is called {\em rigid} if for some
increasing sequence $(n_i)$ of integers, $U_T^{n_i}\to Id$ in
$L^2\xbm$.

\begin{Lemma}\label{elith}
If $(C_2(\mu),\tilde{T})$ is pointwise recurrent then
$(X,\cb,\mu,T)$ is rigid.
\end{Lemma}

\begin{proof}
Let $\alpha=\{A_1,\dots,A_k\}$ be a finite measurable partition of
$X$
such that the numbers $a_i=\mu(A_i)$ are
all positive and distinct ($a_i\neq a_j$ whenever $i\neq
j$). Fix $n \in \N$ and set $t_{ij}=\mu(T^{-n}A_i\cap A_j)$.
Thus $\sum_{i=1}^k a_i=1$
and for each $i$, $\sum_{j=1}^k t_{ij}=a_i>0$.

Let $\xi = \xi_\alpha$ be the measure on $X\times X$ defined by:
$$
\xi = \sum_{i=1}^k a_i \mu_{A_i} \otimes \mu_{A_i}, \quad
{\mbox{where}} \quad \mu_{A_i}(B) = \frac{\mu(B\cap
A_i)}{\mu(A_i)},\;1\leq i\leq k.
$$
If we fix $B\in\cb$ then $\xi(B\times X)=\xi(X \times B)
=\sum_{i=1}^k a_i\frac{\mu(A_i\cap B)}{\mu(A_i)}=\mu(B)$, so
$\xi\in C_2(\mu)$. Notice that $\xi(\bigcup_{j=1}^k A_j\times
A_j)=1$.

Fix an $\vep>0$. Since $\wt{T}$ is pointwise recurrent, we can
assume that for some $n\geq1$,  $\tilde{T}^n\xi$ is so close to
$\xi$ that
$$
|\tilde{T}^n \xi(\cup_{j=1}^k A_j \times A_j) -1|<\vep.
$$
If we now put $b_{ij}=\frac{\mu(T^{-n} A_j\cap A_i)}{a_i}$ then
$\sum_{j=1}^kb_{ij}=1$ and
$$
\begin{array}{l}
\tilde{T}^n \xi(\cup_{j=1}^k A_j \times A_j)\\
= \sum_{j=1}^k \sum_{i=1}^k a_i\mu_{A_i} \otimes \mu_{A_i}
(T^{-n} A_j \times T^{-n} A_j)\\
= \sum_{i=1}^k a_i \sum_{j=1}^k
\left(\frac{\mu(T^{-n} A_j\cap A_i)}{a_i}\right)^2\\
 = \sum_{i=1}^k a_i \sum_{j=1}^k (\frac{t_{ij}}{a_i})^2=
 \sum_{i=1}^k a_i \sum_{j=1}^k b_{ij}^2,
\end{array}
$$
so
$$
\left|\sum_{i=1}^k a_i \sum_{j=1}^kb_{ij}^2 -1\right|<\vep.
$$
If $\vep$ is sufficiently small then this implies that for all
$1\leq i\leq k$, $\left|\sum_{j=1}^kb_{ij}^2 -1\right|<\vep'$ or,
equivalently $\sum_{j\neq j'}b_{ij}b_{ij'}<\vep'$. This, in turn,
means that, given $i$, for only one $j_i$ we have
$|t_{ij_i}-a_i|<\vep^{\prime\prime}$. The map $i\mapsto j_i$ is
1-1. However the numbers $a_i$ are distinct, so if for each $i$
$$
|\mu(T^{-n}A_i\cap A_{ij_i})-\mu(A_i)|<\vep^{\prime\prime}
$$
then $j_i=i$.
It follows that $j_i=i$ for each $1\leq i\leq k$, that is, $|\mu(T^n A_i \cap
A_i)-a_i|<\vep^{\prime\prime}$. This proves the rigidity of
$(X,\cb,\mu,T)$ since $\vep^{\prime\prime}\to 0$ when $\vep\to0$.
\end{proof}

\begin{Th} \label{sztywny}
The following conditions on $T \in \Aut\xbm$ are equivalent:\\
(i) $T$ is rigid;\\
(ii) $\wt{T}$ is pointwise recurrent;\\
(iii) $\wt{T}$ is rigid;\\
(iv)  $\wt{T}$ is uniformly rigid.
\end{Th}

\begin{proof}
(i)$\Rightarrow$(iv) Let $\{f_k\}_{k\geq 0}$ be a
dense set of $L^2$-functions all of norm~1 (see~(\ref{compatible})). Let
$J\in\cj(\mu)$ and fix $\vep>0$. Choose  $N\geq1$  so that
$\sum_{k,l=N}^\infty \frac1{2^{k+l}}<\vep/8$. Since $T$ is rigid
we can find $m\geq1$ such that
$$
\sum_{k,l=0}^{N-1}\frac1{2^{k+l}}(\|U_T^mf_k-f_k\|+\|U_T^mf_l-f_l\|)
<\vep/2.$$ We now have
$$
\sum_{k,l\geq0}\frac1{2^{k+l}}|\langle J\circ
U_T^mf_k,U_T^mf_l\rangle -\langle Jf_k,f_l\rangle|\leq
$$
$$
\sum_{k,l\geq0}\frac1{2^{k+l}}(|\langle J\circ
U_T^mf_k,U_T^mf_l\rangle- \langle Jf_k,U_T^mf_l\rangle|+|\langle
Jf_k,U_T^mf_l\rangle-\langle Jf_k,f_l\rangle|)\leq
$$
$$
\sum_{k,l=0}^{N-1}\frac1{2^{k+l}}(\|U_T^mf_k-f_k\|+\|U_T^mf_l-f_l\|)
+\sum_{k,l\geq N
}\frac1{2^{k+l}}(\|U_T^mf_k-f_k\|+\|U_T^mf_l-f_l\|)$$ and
therefore
$$
\sum_{k,l\geq0}\frac1{2^{k+l}}|\langle \wt{T}^m(J) f_k,f_l\rangle
-\langle Jf_k,f_l\rangle|<\vep/2+\vep/8<\vep.
$$
The uniform rigidity of $\wt{T}$ follows.

The implications (iv) $\Rightarrow$ (iii)  $\Rightarrow$ (ii)
are true for every dynamical system and finally the implication
(ii) $\Rightarrow$ (i) follows from Lemma~\ref{elith}.
\end{proof}

\subsection{The discrete spectrum case and distality of the topological
lens}\label{3.5}

Suppose $T$ is ergodic and has discrete spectrum. In
this case the set $\{U_T^n:\:n\in\Z\}$ is relatively compact in
the strong operator topology, and moreover each limit point 
of this set is of the
form $U_S$, where $S\in C(T)$. As we have seen in the proof of
Theorem~\ref{sztywny} this implies that whenever $J\in
\cj(\mu)$ and $U^{n_k}_T\to U_S$ then $\wt{T}^{n_k}(J)\to
U_S^{-1}\circ J\circ U_S$. It follows that any pointwise limit
$\Theta$ of powers of $\wt{T}$, that is any element of the Ellis
semigroup of $\wt{T}$,  is also a conjugation. Hence the Ellis
semigroup of $\wt{T}$ is a group of homeomorphisms. This fact implies
that $\wt{T}$ is equicontinuous, that is, the
family $\{\wt{T}^n:\:n\in\Z\}$ is equicontinuous. Moreover each
minimal subset (which must be the closure of an orbit) is of the form
$$\{U_S\circ J\circ U_{S^{-1}}:\: S\in C(T)\}$$
for some $J\in \cj(\mu)$.

Equicontinuous systems are special examples of distal systems.
Recall that a homeomorphism $R$ of a compact metric space
$Z$ is {\em distal} if for every pair $(z_1,z_2)$ of distinct points of $Z$ the
closure of its orbit (via $R\times R$) is disjoint from the
diagonal $\Delta_Z=\{(z,z):\:z\in Z\}$ in $Z\times Z$. Every distal system
has a decomposition into minimal components, in
other words every point $z\in Z$ is almost periodic (uniformly
recurrent) that is, the set of return times of $z$ to
any fixed neighborhood has bounded gaps.
(For more details see e.g. \cite{Gl}.)

\begin{Th}\label{eli3}
Suppose that $T$ is  ergodic and that $\wt{T}$ on $C_2(\mu)$ is distal. 
Then $T$ has discrete spectrum.
Therefore if $\wt{T}$ is distal then it is equicontinuous.
\end{Th}

\begin{proof}
Using the notation of Lemma~\ref{elith}, by our
assumption of distality,  we obtain that for a finite measurable
partition $\alpha=\{A_1,\dots,A_k\}$ of $X$ with the property that
the $k$ numbers $a_i=\mu(A_i)$ are distinct and positive, the
corresponding measure,
$$
\xi_\alpha = \sum_{i=1}^k a_i \mu_{A_i} \times \mu_{A_i}, \quad
{\mbox{where}} \quad \mu_A(B) = \frac{\mu(B\cap A)}{\mu(A)}
$$
is a uniformly recurrent point of the topological system
$(C_2(\mu),\tilde{T})$. However, as we have seen in the proof of
that lemma, this implies that the set of recurrence times for the
sets $A_i$ ($i=1,2,\dots,k$), i.e. the set
$$
\{n\in\Z: |\mu(T^n A_i \cap A_i)-\mu(A_i)|< \vep\},
$$
has bounded gaps or equivalently the set $\{n\in\Z:
\|U_T^{n_i}1_{A_i}-1_{A_i} \|_2<\sqrt2 \vep;\ i=1,2,\dots,k\}$ has
bounded gaps. It follows that for each $\vep>0$ the set $
\{U_T^n(1_A):\:n\in\Z\} $
admits a finite $\vep$-net, which means that its
$L^2$-closure is compact. In turn this implies that for every
$k$-tuple of real numbers $(c_1,c_2,\dots,c_k)$ the
$L^2(\mu)$-function
$$
\sum_{i=1}^k c_i \chi_{A_i}
$$
also has a compact $U_T$-orbit and therefore $T$ has
discrete spectrum (see \cite{Ku} or notice that in the terminology of
\cite{Fu} we obtained a dense set of compact functions  in
$L^2(\mu)$). This completes the proof of the theorem.
\end{proof}

In particular, for $T$ with non-discrete spectrum there
always are non-trivial proximal pairs in the topological lens.

\vspace{0.3cm}

\begin{Prop}\label{DS}
The topological lenses of aperiodic ergodic rotations are topologically conjugate iff the rotations are conjugate as measure preserving systems. 
\end{Prop}

\begin{proof}
We first point out the following facts concerning an aperiodic ergodic rotation
$T \in \Aut(\mu)$.
 
1.\ 
Each minimal subset of $\widetilde{T}$ is of the form
$Y_J=\{U_S\circ J\circ U_{S^{-1}}:\:S\in C(T)\}$ for $J\in {\cal
J}(\mu)$.

2.\ 
Assume additionally (but with no loss in generality) 
that $Tx=x+x_0$, i.e. $T$ is
a uniquely ergodic rotation on a compact monothetic (metric) group $X$.
Then  $T$ is topologically conjugate to the translation  $t_T$ by
$T$ on $C(T)$.

3.\ 
Under the assumption in 2., $\widetilde{T}|_{Y_J}$ is a
topological factor of $T$. Indeed, the map
$$
C(T)\ni S\mapsto U_S\circ J\circ U_{S^{-1}}\in Y_J
$$ 
is equivariant (between $t_T$ and $\widetilde{T}|_{Y_J}$).

4.\ 
If $J$ does not commute with any $S\in C(T)\setminus\{Id\}$ then the 
above map is an isomorphism.

5.\ There always is some $J \in {\cal{J}}(\mu)$ 
which does not commute with any $S\in C(T)\setminus\{Id\}$, e.g. take $J = R$, where $R\in \Aut(\mu)$ is weakly mixing with $C(R) = \{R^n: n \in \Z\}$.

Now assume that for aperiodic ergodic rotations $T_1$ and $T_2$
we have a topological conjugacy $\Phi: (C_2(\mu),\widetilde{T_1}) \to
(C_2(\mu),\widetilde{T_2})$. Then $\Phi$ sends minimal
subsets onto minimal subsets and moreover every minimal subset
has a minimal preimage. It now follows, in view of (1) - (5),
that $T_1$ and $T_2$ are weakly topologically isomorphic (i.e. each is a
topological factor of the other), hence are isomorphic (a well
known fact). Of course this implies that they are also
measure theoretically isomorphic.
\end{proof}

\subsection{Weak mixing of a system versus topological transitivity
of its lens}\label{Sec-WM}

It is well-known that invertible elements are dense
in the set $C_2(\mu)$ (see \cite{Gl-Ki}). However we will need a
slightly more concrete result saying that a special family of
interval exchange transformations is dense in
$C_2(\lambda_{[0,1]})$, where $\la_{[0,1]}$ stands for Lebesgue
measure
on $[0,1]$. 

Thus we now assume
(with no loss in generality)
that $X=[0,1]$, $\mu=\la_{[0,1]}$.
Let $\alpha=\{I_1,I_2\dots, I_k\}$ be the partition of $[0,1]$ into 
$k$ intervals  of equal length. Fix $\vep>0$ and $P$ as in~(\ref{postac}).
We also assume that $p_{ij}\in\Q$ and let
$$
p_{ij}=\frac{m_{ij}}L,\quad m_{ij}\in\N,\quad L\in\N\setminus\{0\}
$$
for all $i,j=1,\ldots,k$. Divide each interval $I_i$ into $L$
subintervals  $J_{ij}$ of equal length. We are now going to define
a transformation $S\in \Aut([0,1],\lambda_{[0,1]})$ which will be an
element of $U(\alpha,\epsilon,P)$. It will be defined as an
interval exchange automorphism; i.e.\ each $J_{ij}$ will be
mapped by $S$ onto an interval $J_{\sigma(i,j)}$ via a
map of the form $x\mapsto x+\beta_{ij}$ (defined on $J_{ij}$),
where $\sigma$ is a suitable bijection of
$\{(i,j):\:1\leq i\leq k,1\leq j\leq L\}$. In fact, we will group
some consecutive $J_{ij}'s$ into longer subintervals and then
permute these new subintervals.
Therefore if we ``visualize" the graph of $S$ in $[0,1]\times[0,1]$ as
given by the diagonals of some little subsquares, we only
need to say what is this family of subsquares.
We begin by taking the subsquares:
$$
(J_{11}\cup\ldots
\cup J_{1,km_{11}})\times (J_{11}\cup\ldots \cup J_{1,km_{11}}),$$$$
(J_{1,km_{11}+1}\cup\ldots \cup J_{1,km_{11}+km_{12}})\times
(J_{21}\cup\ldots \cup J_{2,km_{12}}),\ldots,
$$
$$
(J_{1,km_{11}+km_{12}+\ldots+km_{1,k-1}+1}\cup\ldots \cup
J_{1,km_{11}+km_{12}+\ldots+km_{1k}})\times (J_{k1}\cup\ldots \cup
J_{1,km_{1k}})
$$
which define $S$ on $I_1$ (since $\sum_{j=1}^k
m_{1j}/L=1/k$, $\sum_{j=1}^k km_{1j}=L$). To define $S$ on $I_2$
we choose the following subsquares:
$$
(J_{21}\cup\ldots \cup J_{2,km_{21}})\times
(J_{1,km_{11}+1}\cup\ldots \cup J_{1,km_{11}+km_{21}}),
$$
$$
(J_{2,km_{21}+1}\cup\ldots \cup J_{2,km_{21}+km_{22}})\times
(J_{2,km_{12}+1}\cup\ldots \cup J_{2,km_{12}+km_{22}}),\ldots
$$
$$
\begin{array}{l}
(J_{2,km_{21}+\ldots +km_{2,k-1}+1}\cup\ldots \cup
J_{2,km_{21}+\ldots+km_{2,k-1}+km_{2k}})\times\\
\qquad \qquad (J_{k,km_{1k}+1}\cup\ldots \cup J_{k,km_{1k}+km_{2k}}).
\end{array}
$$
We keep going with this procedure of choosing ``the first possible"
square, defining $S$ on $I_3$, then through all the remaining intervals. This
construction is correct since for each $a,b=1,\ldots,k$
$$
\sum_{j=1}^k\frac{km_{aj}}L=1,\quad \sum_{i=1}^k\frac{km_{ib}}L=1.
$$

Given $k\geq1$, let $\cs_k$ stand for the family of automorphisms
of $([0,1],\la)$ given by dividing $[0,1]$ into $k$ intervals of
equal length and then permuting them according to a permutation
$\pi$ of $\{1,2,\ldots,k\}$. By  the above reasoning we have
proved the following:

\begin{Prop}\label{IET}
The family $\bigcup_{k\geq1}\cs_k$ is
dense in $C_2(\la_{[0,1]})$. \bez
\end{Prop}

\begin{Remark}\em
It has been proved by
Kechris and Rosendal in \cite{KR} that there
exists a residual set of $T$'s, such that $\wt{T}|{\Aut\xbm}$,
that is, conjugation by $T$, is topologically transitive on $\Aut\xbm$.
Hence by the above remark we also have that for a
residual set of $T$'s, $\wt{T}$ is transitive on $C_2(\la_{[0,1]})$.
As the set of weakly mixing transformations is residual in  $\Aut\xbm$,
it follows that for a residual set of weakly mixing
transformations $T$, conjugation by $T$ is transitive. It turns out, however,
that this property is in fact a characterization of weak mixing as the
theorem below shows.
\end{Remark}

\begin{Th}\label{slabemiesz}
Assume that $T \in \Aut\xbm$ is ergodic.
Then the following conditions are equivalent.\\
(i) $T$ is weakly mixing.\\
(ii) $\wt{T}$ is transitive.\\
(iii) $\wt{T}$ is topologically weakly mixing.
\end{Th}

\begin{proof}
(ii)$\Rightarrow$(i) 
Denote by $(X_1,\mu_1,T_1)$ the Kronecker factor of $(X,\mu,T)$. 
In view of~(\ref{w2}), $\wt{T_1}$ is a factor of $\wt{T}$ and
since by assumption $\wt{T}$ is transitive, so is
$\wt{T_1}$. Since $\wt{T_1}$ is WAP it has exactly one
minimal set. However, if it is not the trivial one point system,
it has at least two distinct fixed points, $\mu_1 \times \mu_1$ and
$\Delta_{\mu_1}$. Thus $\wt{T}_1$ is trivial,
hence $T$ is weakly mixing.

(i)$\Rightarrow$(ii) For a partition $\alpha=\{A_1,\ldots,A_k\}$
with $\mu(A_j)=1/k$ for $j=1,\ldots,k$, a permutation $\eta$ of
$\{1,\ldots,k\}$ and an $\vep>0$ set
$$
V(\alpha,\eta,\vep)=\{\xi\in C_2(\mu):\:\left|\xi(A_i\times
A_{\eta(i)})-\frac1k\right|<\vep\}.
$$
{\bf Claim.} Given two
permutations $\pi,\sigma$ of $\{1,\ldots,k\}$ there exists an
$n\geq1$ such that \beq\label{em1}
\wt{T}^{-n}(V(\alpha,\sigma,\vep))\cap
V(\alpha,\pi,2\vep')\neq\emptyset,\eeq where $\vep'>0$ is made
precise below.

Indeed, by the weak mixing property of $T$ there exists $n\geq1$
such that
\beq\label{em2}
\left|\mu(T^nA_i\cap
A_j)-\frac1{k^2}\right|<\vep\;\;\mbox{for
all}\;\;i,j=1,\ldots,k.
\eeq
Let
$$
\Delta_\eta=\bigcup_{i=1}^k A_i\times
A_{\eta(i)},\;\mbox{where}\;\eta=\sigma,\pi.
$$
For a fixed $1\leq i\leq k$ let $\ov{B}^i_s=T^{-n}A_s\cap A_i$,
$s\geq1$. Then $\beta_i=\{\ov{B}^i_1,\ldots,\ov{B}^i_k\}$ is a
partition of $A_i$ and by~(\ref{em2})
$$
\left|\mu(\ov{B}^i_s)-\frac1{k^2}\right|<\vep\;\;\mbox{for}\;s=1,\ldots,k.
$$
We now replace the partitions $\beta_i$ by partitions for which
all atoms $B^i_s$, $1\leq s\leq k$, have measure $\frac1{k^2}$ and
moreover \beq\label{oszac} \sum_{i,s=1}^k\mu(\ov{B}^i_s\triangle
B^i_s)<\vep',\eeq where $\vep'=\vep'(\vep,k)$ and $\vep'\to0$ when
$\vep\to0$. Now choose any $\xi\in C_2(\mu)$ such that
$$
\xi(B^i_s\times  B^{\sigma(i)}_{\pi(s)})=\frac1{k^2}
$$
for all
$i,s=1,\ldots,k$ (for example we can take the measure
$\sum_{i,s=1}^k\frac1{k^2}\mu_{B^i_s}\ot\mu_{B^{\sigma(i)}_{\pi(s)}}$,
see the proof of Lemma~\ref{elith}). Notice that the measure
$\xi$ is supported by the union $\bigcup_{i,s=1}^k B^i_s\times
B^{\sigma(i)}_{\pi(s)}$, and in particular \beq\label{em3}
\xi(B^i_s\times
B^j_t)=0\;\;\mbox{unless}\;j=\sigma(i),t=\pi(s).\eeq
We will now check that\\
(i) $\xi\in V(\alpha,\sigma,\vep)$,\\
(ii) $\wt{T}^n\xi\in V(\alpha,\pi,2\vep')$.

Indeed, as $\bigcup_{s=1}^kB^i_s=A_i$ for each $1\leq i\leq k$,
and~(\ref{em3}) holds,
$$
\xi(A_i\times A_{\sigma(i}))=\xi\left(\bigcup_{s=1}^kB^i_s\times
B^{\sigma(i)}_{\pi(s)}\right)=
$$$$
\sum_{s=1}^k\xi\left(B^i_s\times
B^{\sigma(i)}_{\pi(s)}\right)=k\cdot\frac1{k^2}=\frac1k
$$
and thus
$\xi\in V(\alpha,\sigma,\vep)$. Moreover,
$\xi\left(\bigcup_{i=1}^k A_i\times A_{\sigma(i)}\right)=1$. In
order to check~(ii), consider
$$
(\wt{T}^n\xi)(A_s\times A_t)=
\xi(\bigcup_{i=1}^k(A_i\cap
T^{-n}A_s)\times\bigcup_{j=1}^k(A_j\times T^{-n}A_t))=
$$
$$
\xi(\bigcup_{i,j=1}^k\ov{B}_{s}^i\times \ov{B}^j_t)
=\sum_{i,j=1}^k \xi(\ov{B}_{s}^i\times \ov{B}^j_t).
$$
We have
$$
\left| \sum_{i,j=1}^k \xi(\ov{B}_{s}^i\times \ov{B}^j_t)-
\sum_{i,j=1}^k \xi(B_{s}^i\times B^j_t)\right|\leq
$$
$$
\sum_{i,j=1}^k |\xi(\ov{B}_{s}^i\times \ov{B}^j_t) -
\xi(B_{s}^i\times B^j_t)|\leq \sum_{i,j=1}^k
(\mu(\ov{B}_{s}^i\triangle B^i_s)+\mu(\ov{B}^j_t\triangle
B^j_t))<2\vep'.
$$
In particular,
$$
\left|(\wt{T}^n\xi)(A_s\times
A_{\pi(s)})-\sum_{i,j=1}^k\xi(B^i_s\times
B^j_{\pi(s)})\right|<2\vep'.
$$
But in view of (\ref{em3})
$$
\sum_{i,j=1}^k\xi(B^i_s\times B^j_{\pi(s)})=
\sum_{i=1}^k\xi(B^i_s\times B^{\sigma(i)}_{\pi(s)})=\frac1k
$$
and the proof of the claim is complete.

Let $U_1,U_2$ be any non-empty open subsets of $C_2(\mu)$. By
Proposition~\ref{IET} there are a partition
$\alpha=\{A_1,\ldots,A_k\}$ of $X$ with $\mu(A_i)=1/k$,
$i=1,\ldots, k$, $\vep>0$ and permutations $\pi,\sigma$ such that
$$
V(\alpha,\sigma,\vep)\subset U_1,\;\;V(\alpha,\pi,2\vep')\subset U_2.
$$
By the above claim there exists $n\geq1$ such that
$$
\wt{T}^{-n}V(\alpha,\sigma,\vep)\cap
V(\alpha,\pi,2\vep')\neq\emptyset
$$
and therefore
$\wt{T}^{-n}U_1\cap U_2\neq\emptyset$ which completes the proof of
this part of the theorem.

(i)$\;\Rightarrow\;$(iii) In order to obtain topological weak
mixing we repeat the arguments used for the proof of transitivity but
now with
$$
\wt{T}^{-n}V(\alpha,\sigma_i,\vep)\cap
V(\alpha,\pi_i,2\vep')\neq\emptyset
$$
for $i=1,2$.
\end{proof}

\begin{Remark}\em
Recall that the map $S\mapsto \mu_S$ is a
homeomorphism from the Polish group $\Aut\xbm$ onto a dense
$G_\delta$ subset of $C_2(\mu)$ which is $\wt{T}$-invariant.
Moreover it intertwines $\ov{T}$ and $\wt{T}$, where by $\ov{T}$ we
denote the action of $T$ on $\Aut\xbm$ by conjugation (see \cite{Gl-Ki}). 
We conclude that $T$ is weakly mixing iff $\ov{T}$ is topologically transitive
(iff it is topologically weakly mixing) on the Polish space $\Aut\xbm$.
\end{Remark}

For $T,S\in \Aut\xbm$ denote by $\langle T,S\rangle$ the closed
subgroup generated by $T$ and $S$.

\begin{Cor}\label{Eli2}
If $T$ is weakly mixing then
$$
\mathcal{S}=\{S\in \Aut\xbm:\:\langle T,S\rangle=\Aut\xbm\}
$$
is a residual subset of $\Aut\xbm$.
\end{Cor}

\begin{proof} 
By Theorem \ref{slabemiesz} 
$$
\{S\in \Aut\xbm:\:\{T^{-n}ST^n:\:n\in\Z\}\;\mbox{is dense
in}\; \Aut\xbm\}
$$
is dense $G_\delta$ and it is clearly a subset of $\mathcal{S}$.
\end{proof}

\begin{Remark}\em
The proof of Theorem~\ref{slabemiesz} shows also that if
$T$ is mixing then $\wt{T}$ is topologically mixing.
\end{Remark}

\section{Periodic points of the lens}

Recall that a topologically transitive dynamical system $(Y,S)$ is
a {\it $P$-system} (or is {\it chaotic in the sense of Devaney})
if the set of $S$-periodic points is dense in $Y$
(see \cite{Gl-We-93}).

Let  $T\in \Aut\xbm$. We recall that an element $\xi \in
C_2(\mu)$ is a fixed point for $\tilde{T}$ ($\tilde{T}\xi = \xi$)
if and only if $\xi$ is a self-joining for $T$, i.e. $\xi\in
J_2(T)$. Thus $\xi \in C_2(\mu)$ is a periodic point for
$\tilde{T}$ ($\tilde{T}^n\xi = \xi$ for some $n \ge 1$) if and
only if $\xi\in J_2(T^n)$. Set
\begin{gather*}
J^n(T)=\{\xi \in C_2(\mu):
\tilde{T}^n\xi = \xi\},  \\
J^\infty(T) = \bigcup \{J^n(T): n \ge 1\}.
\end{gather*}
Similarly we let
\begin{gather*}
C^n(T)=\{S \in \Aut(X,\cb,\mu): S T^n=T^n S\},   \\
C^\infty(T) = \bigcup \{C^n(T): n \ge 1\}.
\end{gather*}
(Thus $C^1(T)=C(T)$ is the centralizer of $T$ in $\Aut(X,\cb,\mu)$,
and for each $n\ge 1$, $C^n(T)=C(T^n)$.)
Identifying
a transformation $S \in \Aut(X,\cb,\mu)$ with its graph measure
$\Delta_S$ (the image of $\mu$ under the map $x \mapsto (x,Sx)$ of
$X$ into $X\times X$), we can think of $\Aut(X,\cb,\mu)$ as a dense
$G_\delta$ subset of $C_2(\mu)$. In particular, viewed in this
way, $C^n(T)$ is a subset of $J^n(T)$ for $n=1,2,\dots$.


\subsection{The topological lens of a $K$-transformation is
a $P$-system}

\begin{Th}\label{elibenji1}
Let $(X,\cb,\mu,T)$ be an ergodic system.
\begin{enumerate}
\item
If $(X,\cb,\mu,T)$ is Bernoulli then the set $C^\infty(T)$ is
dense in $\Aut(X,\cb,\mu)$.
\item
If $(X,\cb,\mu,T)$ is a $K$-system then the set $J^\infty(T)$ is
dense in $C_2(\mu)$. Hence, if $(X,\cb,\mu,T)$ is a $K$-system
then the topological system $(C_2(\mu),\tilde{T})$ is a
$P$-system.
\end{enumerate}
\end{Th}

\begin{proof}
Let $T$ be a Bernoulli transformation. Given a measurable
partition of $X$, $\ca=\{A_0,A_1,\dots,A_{d-1}\}$ into $d$ sets of
equal measure, and an $\vep>0$, we would like to find an $n_0$ and
$S \in C(T^{n_0})=C^{n_0}(T)$ such that
\begin{equation}
\mu(SA_i\triangle A_{i+1}) < \vep, \quad \mbox{\ for\ }
i=0,\dots,d-1\ \pmod d.\label{S}
\end{equation}


{\bf Step 1:}
As $T$ is $K$, there exists $n_0$ so that the
partitions $\{T^{jn_0} \ca\}_{j\in \Z}$ are $\delta$-independent.
Here, $\delta$ is chosen small enough so that by the Ornstein
version of Sinai's theorem (see \cite{OS}, Lemma 5), there exists a measurable
partition $\hat\ca=\{\hat A_0, \hat A_1,\dots,\hat A_{d-1}\}$ of
$X$ into $d$  sets with $\mu(\hat A_i)=\frac1d$, such that
\begin{equation}
\begin{array}{l}
{\mbox{the partitions}}\ \{T^{jn_0}\hat\ca\}_{j\in \Z}
\ {\mbox{are independent}},\\
\sum_{i=0}^{d-1} \mu(A_i\triangle \hat A_i) <
\frac{\vep}{100}.
\end{array}
\end{equation}

{\bf Step 2:}
Next use Thouvenot's
relative version of Sinai's theorem
(see \cite{Th} and \cite{Kieff}), to get a complementary partition
$\cp=\{P_0,P_1,\dots,P_{\ell-1}\}$ such that
\begin{equation}
\label{hat}
\begin{array}{l}
{\mbox{the partitions}}\ \{T^{jn_0} \cp\}_{j\in \Z} \ {\mbox{
are independent}},\\
\vee_{-\infty}^{\infty} T^{jn_0} \hat \ca \perp
\vee_{-\infty}^{\infty} T^{jn_0} \cp,\\
 H(\cp) + H(\hat \ca) = h(X, T^{n_0}).
\end{array}
\end{equation}

{\bf Step 3:}
Of course $T^{n_0}$ is also Bernoulli and Ornstein's
isomorphism theorem
(see \cite{O}, Proposition 11, page 31)
says that the partition $\cp\vee\hat\ca$ can
be modified by an arbitrarily small amount so as to produce
partitions $\tilde{\cp}$ and $\tilde{\ca}$ with the properties
\begin{equation}
\label{tilde}
\begin{array}{l} {\mbox{the partitions}}\
\{T^{jn_0}(\tilde{\cp}\vee\tilde{\ca})\}_{j\in \Z}
\ {\mbox{are independent}},\\
 \vee_{-\infty}^{\infty} T^{jn_0} \tilde{ \ca} \perp
\vee_{-\infty}^{\infty} T^{jn_0} \tilde{\cp},\\
 \vee_{-\infty}^{\infty}T^{jn_0}(\tilde{\cp}\vee\tilde{\ca} )
\ {\mbox{is the full Borel $\sigma$-algebra of $X$}},\\
 \sum_{i=0}^{d-1} \mu(\tilde{A}_i\triangle \hat A_i) < \frac{\vep}{100}.
\end{array}
\end{equation}

{\bf Step 4:} Now we can view $T^{n_0}$ as a shift on $d+\ell$
symbols $\{\alpha_i: i \in \{0,\dots,d-1\}\}$ and $\{\beta_i: i
\in\{0,\dots,\ell-1\}\}$, with the product measure. On the phase
space we define a transformation $S$ by the map
$$
\left(
\begin{array}{llll} \dots ,&\alpha_{-1},&\alpha_0,&\alpha_1,\dots \\
\dots ,&\beta_{-1},&\beta_0,&\beta_1,\dots
\end{array}
\right) \mapsto \left(
\begin{array}{llll} \dots ,&\alpha_{-1}+1,&\alpha_0+1,&\alpha_1+1,\dots \\
\dots ,&\beta_{-1},&\beta_0,&\beta_1,\dots
\end{array}
\right) \ \pmod d.
$$
This map is measure preserving, commutes with the shift, and has
the property
\begin{equation}
S \tilde{A_i} = \tilde{A}_{i+1} \quad \mbox{\ for\ }
i=0,\dots,d-1\ \pmod d.\label{=}
\end{equation}
Taking into account the formula (\ref{hat}), (\ref{tilde}) and
(\ref{=}), we obtain (\ref{S}). This proves part 1 of the theorem.

To prove part 2 it suffices to show that, given a
$K$-transformation $T$, a partition $\ca=\{A_0, A_1,\dots,
A_{d-1}\}$ of $X$ into $d$ sets with $\mu( A_i)=\frac1d$, and
$\vep>0$, there are $n_0$ and a self-joining $\la$ of $T^{n_0}$
such that
\begin{equation}
\la (\cup_{i=0}^{d-1} A_i \times A_{i+1}) > 1 -\vep \ \pmod d.
\label{K}
\end{equation}

{\bf Step 5:} As in Step 1 above, there exists $n_0$ so that the
partitions $\{T^{jn_0} \ca\}_{j\in \Z}$ are $\delta$-independent,
and an application of Ornstein's theorem yields a measurable
partition $\hat\ca=\{\hat A_0,\hat A_1,\dots,\hat A_{d-1}\}$ of
$X$ into $d$-sets with $\mu(\hat A_i)=\frac1d$, such that
$$
\begin{array}{l}
 \mbox{the partitions}\ \{T^{jn_0}\hat\ca\}_{j\in \Z}
\ {\mbox{are independent}},\\
\sum_{i=0}^{d-1} \mu(A_i\triangle \hat A_i) < \frac{\vep}{100}.
\end{array}
$$

{\bf Step 6:} Now for the Bernoulli $d$-shift corresponding to
$T^{n_0}$ and the independent partition $\hat\ca$ there is an
automorphism $\hat S$ cyclically permuting the sets $\{\hat
A_i:i=1,2,\dots,d-1\}$. This defines a graph joining $\Delta_{\hat
S}$ on the Bernoulli factor defined by $\vee_{-\infty}^{\infty}
T^{jn_0} \hat{\ca}$. The joining $\Delta_{\hat S}$ gives measure
close to 1 to the set $\cup_{i=0}^{d-1} A_i \times A_{i+1}$ since
it gives measure 1 to the set $\cup_{i=0}^{d-1} \hat A_i \times
\hat A_{i+1}$. Finally lift $\Delta_{\hat S}$ to a self-joining
$\la$ of $T^{n_0}$ on $X \times X$ to get (\ref{K}).
\end{proof}

\subsection{Zero entropy and \protect$P$-systems}

We have been unable to find $T$ with zero entropy such that
$(\widetilde{T},C_2(\mu))$ is a $P$-system (in general, there are
$P$-systems with zero topological entropy, see \cite{Gl-We-93}).
A clear case where the periodic points of the lens fail to be dense
is when the equality
\beq\label{psystem}
J^\infty(T)=J^1(T)=J_2(T)
\eeq
holds, as in this case the set of periodic points is closed and nowhere-dense
(see the remarks following equation (\ref{w3}) in Section \ref{Sec-lens}).  
It has been shown in \cite{Ju-Ru} (Theorem 6.1) that the
equality~(\ref{psystem}) holds for all weakly mixing 2-fold simple systems.
(Although Theorem 6.1 in \cite{Ju-Ru} states merely that
$C(T^k)=C(T)$ for all $k \ne 0$,
the proof actually yields the stronger result $J^e(T^k)=C(T^k)=C(T)=J^e(T)$.
In particular (\ref{psystem}) holds.)

Notice also that if $S$ is a root of $T$ satisfying~(\ref{psystem}) then
the set of $S$-periodic points is non-dense as well.
Indeed, assume that $S^k=T$ for some
$k\geq 2$. Take $r\geq1$ and $J\in\cj(\mu)$ so that $J\circ
U_{S^r}=U_{S^r}\circ J$. Then $U_{S^{kr}}\circ J=J\circ
U_{S^{kr}}$ or equivalently $ U_{T^{r}}\circ J=J\circ U_{T^{r}}$
and therefore $J^\infty (S)\subset J_2(T)$.

The equality~(\ref{psystem}) is also satisfied when
\beq\label{psystem1}
T^n\;\;\mbox{has simple spectrum for
each}\;\;n\geq1.
\eeq 
In fact, if a unitary operator $V\in U(H)$ of
a separable Hilbert space $H$  has
simple spectrum then each
bounded linear operator $W\in L(H)$ commuting with $V$ is a (weak)
limit of polynomials in $V$ and therefore the semigroup of such $W$'s
is commutative. It follows that when $S\in \Aut\xbm$ has
simple spectrum the semigroup $\{J_\la:\:\la\in J_2(S)\}$ is also
commutative. We clearly have $J_2(T)\subset J_2(T^n)$
for each $n \ge1$. Assume now that $J\in \cj(\mu)$ and for
some $n$,  $J\circ U_{T^n}=U_{T^n}\circ J$. Since
$U_T\circ U_{T^n}=U_{T^n}\circ U_T$ and $U_{T^n}$ has simple
spectrum, $J\circ U_T=U_T\circ J$, and therefore
$J_2(T)=J_2(T^n)$, so~(\ref{psystem}) follows. 
Notice that the
property~(\ref{psystem1}) is closed under taking both powers and roots
(indeed, for the latter just observe that whenever
$S^{km}$ has simple spectrum then both $S^k$ and $S^m$ have
simple spectra).

\begin{Th}\label{BEM1}
\begin{enumerate}
\item
The generic transformation $T\in \Aut\xbm$ has the property that
for every  $n \in \Z \setminus \{0\}$,\ $T^n$  is of rank one.
Thus the generic $T$ satisfies  condition~(\ref{psystem1})
and in particular its  set of periodic points is closed and nowhere-dense.
Thus the topological lens is not a $P$-system.
\item
The equality~(\ref{psystem}) holds for all weakly mixing 2-fold simple systems.
In particular, for such systems the lens is never a $P$-system
\item
If $T$ is a Gaussian system with simple spectrum then 
it satisfies  condition~(\ref{psystem1}); in particular $\wt{T}$ is
not a $P$-system.
\end{enumerate}
\end{Th}

\begin{proof}
1. It is well known, and not hard to check, that
the set $\mathcal{R}$ of rank one transformations is a dense
$G_\delta$ subset of the Polish group $G=\Aut(\mu)$.
It is also well known that a rank one transformation has simple spectrum
(see e.g. \cite{Gl}, Theorem 16.5).
For each integer $k \ge 1$, let $\pi_k : \Aut(\mu) \to \Aut(\mu)$ be the
map $T \mapsto T^k$. Clearly $\pi_k$ is a continuous map and
it follows that the set $\mathcal{S}_k:= \pi^{-1}_k(\mathcal{R})
=\{T \in G: T^k \in \mathcal{R}\}$ is a $G_\delta$ subset of $G$.
Clearly $\mathcal{R}$, and therefore also $\mathcal{S}_k$, are conjugation-invariant and, as $\mathcal{S}_k$ contains ergodic transformations
(e.g. an irrational rotation), it follows, by Halmos' theorem, that $\mathcal{S}_k$
is a dense $G_\delta$ subset of $G$.
In fact, it can be shown that $\mathcal{S}_k \subset \mathcal{R}$,
(if $T^k$ is rank 1 then so is $T$), but in any case we can take
$\mathcal{R}_\infty =\mathcal{R}\cap \bigcap_{k \ge 1}\mathcal{S}_k$
as the required generic (in fact, dense $G_\delta$) set of transformations,
since $T\in \mathcal{R}$ iff $T^{-1}\in \mathcal{R}$.

2. This has been shown in \cite{Ju-Ru}.

3. Since $T^n$ is a (generalized) Gaussian system whose maximal
spectral multiplicity is bounded by~$n$, it follows that $T^n$ has
simple spectrum. In fact, the maximal spectral multiplicity of a
Gaussian system is either~$1$ or ~$\infty$ (see \cite{Co-Fo-Si}
and \cite{Le-Pa-Th}).
\end{proof}

Of course every Kronecker (that is, ergodic discrete spectrum) system
is simple, but in this class we encounter both types of behavior.

For an irrational rotation $R_\alpha$ on the circle $\T$, we have
$C(R_\alpha)=C(R_\alpha^n)= C^n(R_\alpha)=\{R_\beta: \beta \in
\T\}$ for all $n\not=0$. Thus, by simplicity, also $J_1(R_\alpha)
=J_n(R_\alpha)$ for all $n \not=0$. Of course, $\widetilde{T}$ is
not a $P$-system because of
the absence of transitivity.

The situation changes when $T$ is not totally ergodic.
Considering the dyadic adding machine $(X,\mu,T)$, with
$X=\{0,1\}^\N$ and $T x = x +{\bf{1}}$, where ${\bf{1}}=(1,0,0,\dots)$, one
can easily check that, given any permutation $\pi$ of
$\Z_{2^n}=\{0,1,\dots,2^n-1\}$, the corresponding homeomorphism
$S_\pi$ defined on $X$ by permuting the first $2^n$ coordinates,
commutes with $T^n$ which do not affect the first $2^n$
coordinates of elements of $X$. It follows that the set
$C^\infty(T)$ is dense in $\Aut(X,\cb,\mu)$; i.e. the periodic
points are dense in $(C_2(\mu),\tilde{T})$. Note however that
$(C_2(\mu),\tilde{T})$ is not topologically transitive and
therefore not a $P$-system.

\section{Invariant measures on
the topological lens and quasi-factors}

In this section we will
examine the relation between
$\wt{T}$-invariant probability measures on the topological
lens and quasi-factors. For completeness
we will repeat some arguments from \cite{Gl} concerning the
basic connection between quasi-factors and self-joinings -- this
has to be checked because we must drop the assumption,
made in \cite{Gl}, of ergodicity of the barycenter. Given
a standard Borel space $(Y,\cc)$ we denote by $M(Y)$ the set of
probability measures on it.
Recall that $M(Y)$ carries a natural
Borel structure which is generated by the evaluation maps
$\phi_A: \mu \mapsto \mu(A), \ A \in \cc$.

Assume now that $T$ is an automorphism  (not necessarily ergodic)
of a standard probability Borel space  $\xbm$. Let $P\in
QF(T,\mu)$ be a {\em quasi-factor} of the system $(T,X,\cb,\mu)$, 
i.e.\ $P\in M(M(X))$ is $T$-invariant for the natural induced action of $T$ on measures, and its barycenter $\int_{M(X)}\nu\,dP(\nu)$ equals $\mu$. 
Consider also
$M_{sym}(X^{\Z})$, the set of probability measures on $X^{\Z}$
which are invariant under the group $\cs$ of all finite
permutations. Let $\phi:M(M(X))\to M_{sym}(X^{\Z})$
be defined by the
formula:
\beq\label{deFHS}
\phi(Q)=\ov{Q}^{\infty}:=\int_{M(X)}(\ldots\otimes\theta\otimes\theta
\otimes\ldots)\,dQ(\theta).
\eeq
By the de Finetti-Hewitt-Savage
theorem $\phi$ is an affine isomorphism of these two convex sets
(simplices).
Now consider $J_{sym,\infty}(T,\mu)$, the set of
infinite symmetric self-joinings of $(T,X,\cb,\mu)$ (it is a
closed subset of $M_{sym}(X^{\Z})$ but since $\mu$ is not
necessarily ergodic, in general this set is not a simplex). We
have
\begin{equation}\label{QF}
\phi(QF(T,\mu))=J_{sym,\infty}(T,\mu).
\end{equation}
In fact, if $P\in QF(T,\mu)$, then clearly $\phi(P)\in M_{sym}(X^{\Z})$ and
the fact that its 1-dimensional marginals equal $\mu$ follows
from the barycenter condition on $P$.
In order to see that $\phi(P)$ is $T^\infty$-invariant note that
\begin{align*}
\phi(P)\circ
T^\infty&=\int_{M(X)}(\ldots\otimes\theta\otimes\theta\otimes\ldots)\,
dP\circ T(\theta)\\
& = \int_{M(X)}(\ldots\otimes\theta\otimes\theta\otimes\ldots)\,
dP(\theta) = \phi(P).
\end{align*}
Finally, if $\rho\in J_{sym,\infty}(T,\mu)$ then, by the de Finetti-Hewitt-Savage
theorem, its $\cs$-ergodic decomposition is of the form
$$
\rho=\int_{M(X)}(\ldots\otimes\theta\otimes\theta\otimes\ldots)\,
dQ(\theta)
$$
Clearly $Q\in QF(T,\mu)$ and $\phi(Q)=\rho$. Thus
$\phi$ is onto and (\ref{QF}) is established.
Note that this shows that~(\ref{deFHS}) is in fact the $\cs$-ergodic
decomposition of $\ov{P}^\infty$.

\begin{Prop} \label{Eli}
The quasi-factor $(T,M(X),P)$ is isomorphic
(as a dynamical system) to the factor of
$(T^\infty,X^{\Z},\ov{P}^\infty)$ given by the $\sigma$-algebra of
symmetric sets (i.e. $\cs$-invariant sets).
\end{Prop}

\begin{proof}
Consider the map
$\Phi:(T^\infty,X^{\Z},\ov{P}^\infty)\to (T,M(X),P)$ given in the
following way. Take first the $\cs$-ergodic decomposition of
$\ov{P}^\infty$:
$$
\ov{P}^{\infty}=\int_{M(X)}(\ldots\otimes\theta\otimes\theta
\otimes\ldots)\,dP(\theta).
$$
Now to $\ov{x}\in X^{\Z}$ which
belongs to exactly one ergodic component (the support of some
$\ldots\otimes\theta\otimes\ldots$) we associate the measure
$\theta=\theta_{\ov{x}}$. All we need to show is that this map
establishes a homomorphism. (Since it is
constant on atoms of the
partition given by the $\cs$-ergodic decomposition, the preimage
of $\cb(M(X))$ will be exactly the $\sigma$-algebra of symmetric
sets.)
Given a Borel subset $\Lambda\subset M(X)$ we have
$$
(\ldots\otimes\theta\otimes\theta\otimes\ldots)(\Phi^{-1}(\Lambda))
=
\begin{cases}
1, & {\text{when $\theta=\theta_{\ov{x}}$ and $\ov{x}\in\Lambda$}}\\
0, & {\text{otherwise}},
\end{cases}
$$
and therefore
$$
\ov{P}^\infty(\Phi^{-1}(\Lambda))=P(\Lambda).
$$
Finally, the map $\Phi$ is clearly equivariant.
\end{proof}

We will now discuss the ergodic case.

\begin{Cor}
Assume that $\mu$ and $P\in QF(T,\mu)$ are
ergodic. Then there exists $\rho\in J_{\infty}^e(T,\mu)$ such that
the factor given by the $\sigma$-algebra of symmetric sets of the
system $(T^\infty,X^{\Z},\rho)$ is isomorphic to the dynamical system
$(T,M(X),P)$.
\end{Cor}

\begin{proof}
By our assumption the $\sigma$-algebra of symmetric sets, which (as
a factor of $(T^\infty,X^{\Z},\ov{P}^\infty)$) is isomorphic to
$(T,M(X),P)$, is ergodic. It follows that the same $\sigma$-algebra
is a factor of a.e. $T^\infty$-ergodic component of
$\ov{P}^\infty$. Since $\mu$ itself is ergodic,
$J_{\infty}(T,\mu)$ is a simplex, so the ergodic components of
$\ov{P}^\infty$ are self-joinings. Any such ergodic self-joining
will serve as $\rho$.
\end{proof}

\begin{Remark}\em
The family of dynamical systems given by ergodic
quasi-factors of an ergodic automorphism $T$ is hence contained in
$$
\{(T^\infty,X/\cb^{\otimes\infty}_{sym},\cb^{\otimes\infty}_{sym},
\la):\:\la\in J^e_\infty(T,\mu)\}.
$$
\end{Remark}

\vspace{0.3cm}

Let us come back to the topological lens $\wt{T}$ of an automorphism
$T$.  We will relate $\wt{T}$-invariant measures to quasi-factors
of $T$.  Notice that $C_2(\mu)\subset M(X\times X)$ is a closed
(compact) subset in the weak-$\ast$ topology, hence each $T\times
T$-invariant Borel measure on $C_2(\mu)$ is also a  $T\times
T$-invariant Borel measure on $M(X\times X)$. Let $P\in
M(C_2(\mu))$ be $\wt{T}$-invariant. If we let $\la=\int_{M(X\times
X)}\nu\,dP(\nu)$ be the barycenter of $P$,
then $\lambda\in C_2(\mu)$
and $(C_2(\mu), P,\wt{T})$ is a quasi-factor of $(X \times X, \la,T\times T)$.
Notice that even if $P$ is ergodic, $\la$ need not be
ergodic. (As an easy example consider $P$ being the Dirac measure at a
non-ergodic self-joining of $T$.) However, as we have shown above, the
arguments from \cite{Gl} about the characterization of quasi-factors
as factors of infinite symmetric self-joinings go through, and
therefore the dynamical system $(C_2(\mu), P,\wt{T})$ is a factor of an
infinite self-joining of the system $(X\times X,\la,T\times T)$ (moreover
this infinite self-joining is symmetric, that is, it is an
invariant measure for the group of finite permutations of
$(X \times X)^\Z$).
It follows directly that $(\wt{T},P)$ is a factor of an infinite
self-joining of $(X,\mu,T)$ and since $(C_2(\mu),P,\wt{T})$ is ergodic, it is
also a factor of an infinite ergodic self-joining of $(X,\mu,T)$.

\begin{Cor}\label{wniosek10}
Assume that $(T,X,\cb,\mu)$ is ergodic and let $P\in
M_{\wt{T}}(C_2(\mu))$ be ergodic. Let $\la\in
J_2(T,\mu)$ denote the barycenter of $P$.
Then $P \in QF(T\times T,\la)$ and there exists
$\rho\in J^e_\infty(T,\mu)$ such
that $(C_2(\mu),P,\wt{T})$ is isomorphic to the factor 
of $(X^\Z,\rho,T^\infty)$ given by the  $\sigma$-algebra of sets
invariant under finite permutations of consecutive pairs of coordinates.
If additionally $\la$ is ergodic, then $\rho\in J^e_\infty(T\times
T,\la)$ (the identification is given by grouping every two
consecutive coordinates).
\end{Cor}

\begin{proof}
Apply Proposition~\ref{Eli} and notice that
$\ov{P}^\infty$ can also be seen as an element of
$J_\infty(T,\mu)$ and that $\mu$ is ergodic.
\end{proof}

\subsection{The entropy of \protect$\wt{T}$}
Let us recall that the
class of zero entropy systems is closed under taking joinings,
since each marginal factor is included in the Pinsker factor of
the system given by the joining. Clearly this class is also closed
under taking factors. Therefore using Corollary~\ref{wniosek10}
and the variational principle we obtain the following.

\begin{Th}\label{qf}
If $T$ has zero entropy, then the topological entropy of $\wt{T}$
is zero.\bez
\end{Th}

We will now show that also the converse is true. In fact a
stronger result holds.

\begin{Th}\label{egmlbw}
If $h(T)>0$ then $h_{top}(\wt{T})=\infty$.
\end{Th}

\begin{proof}
{\bf Step 1.} Assume first that $T$ acting on $\xbm$ is Bernoulli
$(1/2,1/2)$ with an independent generator $(A,A^c)$, i.e.\
$\mu(A)=1/2$ and the family
$\{T^k A,T^k{A^c}:\:k\in\Z\}$ is
independent. Consider the map $F:C_2(\mu)\to[0,1/2]^{\N}$ given by
$$
F(\la)=(\wt{T}^n(\la)(A\times A))_{n\geq0}=(\la(T^{-n}A\times
T^{-n}A))_{n\geq0},
$$
where on the space $[0,1/2]^{\N}$ we
consider the one-sided shift $S$:
$S((x_i)_{i\geq0})=(y_{j})_{j\geq0}$, $y_j=x_{j+1}$.
Clearly, $F$ is an equivariant continuous map. We will show that
$F$ is onto. To this end let us first notice that all we need to
show is that for each $n\geq1$, each block
$(b_0,b_1,\ldots,b_{n-1})\in\{0,1/2\}^n$, there exists $\la\in
C_2(\mu)$ such that $F(\la)[0,n-1]=(b_0,\ldots,b_{n-1})$. In fact,
$\{0,1/2\}^{n}$ is the set of extremal points of the convex set $[0,1/2]^{n}$ and
once we have ``realizations" of the extremal points we use the fact
that the map $\la\mapsto F(\la)[0,\ldots,n-1]$ is affine. The
proof goes by induction on~$n$. For a ``realization" of a block
$B=(b_0)$ of length~1 we take $\la$ as the graph coupling given by
an isomorphism of the space $\xbm$, where the isomorphism is
obtained from two measure preserving maps (isomorphisms)
$$
A_0\mapsto A_{b_0+1},\;A_1\mapsto A_{b_0},
$$
where $A_0 =A$ and $A_1 = A^c$, and addition is mod $2$.
Now take $n\geq1$ and
$B=(b_0,\ldots,b_{n-1})$ and suppose that the graph coupling given
by an isomorphism of the space $\xbm$ obtained from $2^{n-1}$
measure-preserving maps (isomorphisms)
$$
A_{i_0}\cap T^{-1}A_{i_{1}}\cap\ldots\cap
T^{-n+2}A_{i_{n-2}}\mapsto A_{j_0}\cap
T^{-1}A_{j_{1}}\cap\ldots\cap T^{-n+2}A_{j_{n-2}}
$$
 ``realizes" the block $(b_0,\ldots,b_{n-2})$. A coupling $\la$ ``realizing"
the block $B$ is then obtained as the graph coupling given by an
isomorphism of the space $\xbm$ obtained from $2^n$ measure
preserving maps (isomorphisms)
$$
A_{i_0}\cap T^{-1}A_{i_{1}}\cap\ldots\cap T^{-n+2}A_{i_{n-2}}\cap
T^{-n+1}A_{i_{n-1}} \mapsto$$$$ A_{j_0}\cap
T^{-1}A_{j_{1}}\cap\ldots\cap T^{-n+2}A_{j_{n-2}}\cap
T^{-n+1}A_{i_{n-1}+b_{n-1}+1}.
$$
It follows that the shift $S$ is
a topological factor of $\wt{T}$, so $h(\wt{T})=+\infty$.


{\bf Step 2.} Suppose $h(T)>0$. By replacing $T$ with $T^m$ if
necessary, we can assume that the entropy of $T$ is larger than
$\log 2$. Then, by Sinai's theorem, $T$ has a  factor $T_1$ which is Bernoulli
$(1/2,1/2)$. Now $h(\wt{T_1})=+\infty$, so $h(\wt{T})$ is also
infinite.
\end{proof}

\subsection{Invariant measures for $\tilde{T}$
supported on $\Aut(\mu)$}

We start with the following general setup.

\vspace{0.2cm}

{\bf Example:}\
Let $Y$ be a compact abelian second countable topological group with normalized Haar measure $\lambda$. Let $T: Y \to Y$ be a continuous automorphism.
We recall that the measure $\la$ is preserved by $T$; i.e.
$T \in \Aut(Y,\la)$.

With each $z\in Y$ we associate the
$\lambda$-preserving invertible transformation $R_z: Y \to Y$
defined by $R_z(y) = y+ z$. We then have
$T \circ R_z \circ T^{-1}(y) =T \circ R_z ( T^{-1}y)
= T(T^{-1}y + z) = y + Tz =  R_{Tz}(y)$ for all $y, z \in Y$, so that
$$
T \circ R_z \circ T^{-1}=R_{Tz}.
$$
Thus the map $\phi: z \mapsto R_z$ is
a topological isomorphism of the compact
topological dynamical system $(Y,T)$ into the (Polish) dynamical
system $(\Aut(\la),\tilde{T})$, where $\tilde{T}$ denotes conjugation
by $T \in \Aut(\la)$.

Next consider the simplex $M_T(Y)$  of $T$-invariant Borel
probability measures on $Y$. For an element $\nu \in M_T(Y)$
let $P_\nu=\phi_*(\nu)$ be the push-forward image of $\nu$
in the space of probability measures on $\Aut(\la)$ under $\phi$.
It then follows that $\phi$ is an isomorphism of measure dynamical systems
$$
\phi: (Y,\nu,T) \to (\Aut(\la),P_\nu,\tilde{T}).
$$
Identifying $\Aut(\la)$ with its image in $C_2(\mu)$
under the canonical embedding, we obtain the measure
theoretical isomorphism
$$
\phi: (Y,\nu,T) \to (C_2(\mu),P_\nu,\tilde{T}).
$$

Now a particular instance of the above example will establish the
following surprising result.
Let $Y=\T^\Z$, where $\T=\R/ \Z$ is the circle group. Let
$T$ be the shift transformation on $Y$ and let $\la$ be the
product measure $\la =  \la_0^\Z$, where $\la_0$ is normalized
Lebesgue measure on $\T$.
We consider $Y$ also as a compact topological group and
observe that (i) $\la$ is the normalized Haar measure on $Y$ and
(ii) $T$ is a continuous automorphism of the compact group $Y$.

\begin{Th}
The topological lens $(C_2(\la),\tilde{T})$ is universal both topologically
and measure theoretically:
\begin{enumerate}
\item
Every metric compact topologically transitive system appears as a subsystem of
the $G_\delta$ dense $\tilde{T}$-invariant subset $\Aut(\la)\subset C_2(\la)$.
\item
Let $(\Omega,\mathcal{A},\nu,S)$ be an ergodic system. There exists a
$\tilde{T}$-invariant ergodic probability measure $P_\nu$ on $\Aut(\la)$ such that the corresponding dynamical system $(C_2(\mu),P_\nu,\tilde{T})$ is measure theoretically isomorphic to $(\Omega,\mathcal{A},\nu,S)$.
\end{enumerate}
\end{Th}

\begin{proof}
Both claims follow immediately from the example above applied to the
group $Y=\T^\Z$, and the fact that the topological Bernoulli system
$(Y, T)$ is universal both for compact metrizable topologically transitive
systems as well as for ergodic measure preserving systems.
\end{proof}

\begin{Remark}\em
In \cite{Da} Danilenko introduced the notion of near
simplicity. An ergodic dynamical system $(X,\mathcal{B},\mu,T)$ is
{\em 2-fold nearly simple} if each 2-fold ergodic self-joining $\la$ of $\mu$
is either the product measure $\mu\times \mu$ or it is an integral of graph couplings:
$$
\la=\int_{\Aut\xbm}\mu_{S}\,dP(S),
$$
where $P$ is a $\wt{T}$-invariant measure.
For all the examples of 2-fold near-simple maps given in \cite{Da},
the measure $P$ is supported on a finite set.
Danilenko shows that, at least when $P$ is supported on a finite set,
the system $(X \times X, \la, T\times T)$ is
isomorphic to the product of the original system $(X,\mathcal{B},\mu,T)$
and the system $(\Aut\xbm, P, \wt{T})$.
It will be interesting to find near simple systems where $P$ is
a continuous measure.
\end{Remark}

\subsection{$\wt{T}$-invariant measures when \protect$T$
is measure-theoretically distal}
The computation of the simplex $M_{\wt{T}}(C_2(\mu))$ for a
general transformation $T$ seems to be a difficult task.
If $(X,\cb,\mu,T)$ is distal, then by
Corollary~\ref{wniosek10}, $(C_2(\mu),P,\wt{T})$ is
(measure-theoretically) distal for each ergodic $P\in
M_{\wt{T}}(C_2(\mu))$.
In this section we will consider mainly the particularly simple
situation of an isometric extension of a Kronecker system.

Suppose first that $T$ is a Kronecker system, that is, it is ergodic and has
discrete spectrum (and therefore is isomorphic to an ergodic translation on 
a compact metric monothetic group). Assume moreover that 
$P\in M_{\wt{T}}(C_2(\mu))$ is ergodic for $\wt{T}$. 
Then there exists
$J=J_\rho\in\cj(\mu)$ such that $P$ is supported on theÄ set
$$
Y_{J}:=\{U_S\circ J\circ U^{-1}_S:\:S\in C(T)\}
$$
(see the beginning of section~\ref{3.5}). Consider the map
$\kappa: C(T) \to Y_{J}$ defined by \beq\label{FC99} \kappa: S
\mapsto U_S\circ J\circ U^{-1}_S \eeq (notice that it is $1-1$
unless for some $\Id \neq R \in C(T)$ the measure $\rho$ is $R\times
R$-invariant; all measures which are $R\times R$-invariant are
easy to  describe). If on $C(T)$ we consider the action  by the
translation of $T$ then the system we obtain is isomorphic to the
original system $(X,T)$, in particular it is uniquely ergodic and
has discrete spectrum. The system $(Y_{J},\wt{T})$ is a
topological factor of $T$ and hence $(C_2(\mu),P,\wt{T})$ has also
discrete spectrum. In this way we fully described ergodic
invariant measures for Kronecker systems.


Let now $T$ be an arbitrary ergodic automorphism of a standard
probability Borel space $\xbm$. Let $G$ be a compact metric
abelian group and let $\va:X\to G$ be a measurable map (a
cocycle). Define $T_{\va}:X\times G\to X\times G$ by the formula
$$
T_{\va}(x,g)=(Tx,g+\va(x))$$ and observe that $T_{\va}$ preserves
the product measure $\mu\ot\la_G$, where $\la_G$ is the normalized
Haar measure on $G$. We will assume that $T_{\va}$ is ergodic and
has the same eigenfunctions as $T$. (The latter assumption does
not restrict the generality of Proposition~\ref{FC102} below
since, by enlarging the $\sigma$-algebra $\cb$ so that the
additional eigenfunctions become $\cb$-measurable, we can achieve
this situation.) Under these assumptions there exists a
subsequence  $(n_i)$ of density~1 such that $U_{T_{\varphi}}^{n_i}\to 0$
weakly on the space $L^2(X\times G,\mu\ot\la_G)\ominus (L^2\xbm\ot
1_G)$. Given $\chi\in\widehat{G}$ ($\widehat{G}$ stands for the
dual of $G$), we put
$$
V_{T,\va,\chi}(f)=\chi(\va)\cdot f\circ T$$ for each $f\in
L^2\xbm$. Then for each $\chi\neq1$
$$
V_{T,\va,\chi}^{n_i}\to 0\;\;\mbox{weakly on}\; L^2\xbm$$ or
equivalently \beq\label{dz1} \int_X\chi(\va^{(n_i)})\cdot f\circ
T^{n_i}\cdot g\,d\mu\to 0\eeq for each $f,g\in L^2\xbm$.

Suppose now that $\xi\in C_2(\mu\ot\la_G)$ satisfies the condition 
$\xi|_{X\times X}=\Delta_X$. 
Identifying the two copies of $X$ we can assume that
$\xi\in M(X\times G\times G)$ and the  projections of $\xi$ on the
two ``copies'' of $X\times G$ equal $\mu\ot\la_G$. Instead of
$T_{\va}\times T_{\va}$ we must consider $T_{\va\times \va}$. The
disintegration of $\xi$ over $\mu$ has the form $\xi=\int_X
\delta_x\ot\xi^x\,d\mu(x)$, where $\xi^x$ are probability measures
on $G\times G$. Since for each $A\in \cb$, $B\in \cb(G)$
$$\int_A\la_G(B)\,d\mu(x)=\mu(A)\la_G(B)=\xi(A\times B\times G)=\int_{A}\xi^x(B\times G)\,\mu(x),$$
we have $\la_G(B)=\xi^x(B\times G)$. Similarly
$\la_G(B)=\xi^x(G\times B)$, whence $\xi^x\in C_2(\la_G)$ for
$\mu$-a.e.\ $x\in X$.

For $\chi,\eta\in\widehat{G}$ we have
$$
\int_{X\times G\times G}f(x)\chi(g)\eta(h)\,d(T_{\va\times
\va}^{n_i})_\ast(\xi)(x,g,h)=
$$
$$
\int_{X\times G\times
G}f(T^{n_i}x)\chi(g+\va^{(n_i)}(x))\eta(h+\va^{(n_i)}(x))\,d\xi(x,g,h)=
$$
$$
\int_X
(\chi\cdot\eta)(\va^{(n_i)}(x))f(T^{n_i}x)\left(\int_{\{x\}\times
G\times G}\chi(g)\eta(h)\,d\xi^x(g,h)\right)\,d\mu(x).$$ If
$\chi\neq\ov{\eta}$ then, by~(\ref{dz1}), the limit of the latter
expression, when $i\to\infty$, is zero. If $\chi=\ov{\eta}$,
assuming additionally that $T^{n_i}\to S\in C(T)$ and that
$\left(T_{\va\times\va}^{n_i}\right)_\ast(\xi)\to\wt{\xi}$ we
obtain as the limit
$$
\int_{X\times G\times
G}f(x)\chi(g)\eta(h)\,d\wt{\xi}(x,g,h)=\int_{X\times G\times
G}f(Sx)\chi(g)\ov{\chi(h)}\,d\xi(x,g,h).$$ Thus the coupling
$\wt{\xi}$ has the form
$\wt{\xi}=\int_X\delta_{Sx}\ot\wt{\xi}^x\,d\mu(x)$, where for
$\mu$-a.e. $x\in X$ the measure $\wt{\xi}^x$ on $G\times G$ has
the following Fourier coefficients
\beq\label{FC100}\widehat{\wt{\xi}^x}(\chi,\eta)=\left\{\begin{array}{cll}
0, & \mbox{if} & \chi\neq\ov{\eta},\\
\int_{G\times G}\chi(g-h)\,d\xi^{S^{-1}x}(g,h), &\mbox{if} &
\chi=\ov{\eta}.\end{array}\right. \eeq

In order to gain a better picture of the limit measure and its
dynamics we perform the following change of coordinates. Let
$$
m:G\times G\to G,\;\;m(g,h)=g-h\;\;\mbox{and}$$
$$
\ov{m}:X\times G\times G\to X\times G\times
G,\;\;\ov{m}(x,g,h)=(x,g-h,h).$$ Set $\nu^x:=m_\ast(\wt{\xi}^x)$
and $\nu:=\ov{m}_\ast(\wt{\xi})$. We then have 
\beq\label{FC101}
\nu=\int_X\delta_{Sx}\ot\nu^x\ot\la_G\,d\mu(x).\eeq 
In fact,
denoting the measure on the right hand side of~(\ref{FC101}) by
$\xi_1$, we have
$$
\int_{X\times G\times G}f(x)\chi(g)\eta(h)\,d\xi_1(x,g,h)=
$$
$$
\int_Xf(Sx)\left(\int_{\{x\}\times G\times
G}\chi(g)\eta(h)\,d\nu^x(g)d\la_G(h)\right)d\mu(x)=0
$$ 
whenever
$\eta\neq1$ (indeed, $\chi\neq\ov{(\ov{\chi}\cdot\eta)}$ if and
only if $\eta\neq1$). For $\eta=1$ we obtain
$$
\int_{X\times G\times
G}f(x)\chi(g)\eta(h)\,d\xi_1(x,g,h)=\int_Xf(Sx)\left(\int_{\{x\}\times
G} \chi(g)\,d\nu^x(g)\right)d\mu(x)=
$$
$$
\int_Xf(Sx)\left(\int_{\{x\}\times G\times
G}\chi(g-h)\, d\tilde{\xi}^x (g,h)\right)\,d\mu(x)
$$
and we see that indeed
$\ov{m}_\ast(\wt{\xi})=\xi_1$. We also note that
$$
\ov{m}\circ T_{\va\times\va}=T_{0\times\va}\circ\ov{m},
$$ 
where $T_0=T\times Id_G$.

From now on we assume that $T$ is
an ergodic rotation (so
that $C(T)$ is a compact group in the strong operator topology).
We have proved that in this case, under the change of coordinates, the
limit points along subsequences of $(n_i)$ are measures of the
form
$$\int_X\delta_{Sx}\ot\nu^x\ot\la_G\,d\mu(x)=\theta\ot\la_G,$$
where $\theta\in M_\mu(X\times G)$ and the latter is the subset of
$M(X\times G)$ consisting of measures whose projection on $X$ is
$\mu$. In the new coordinate system the automorphism
$T_{\va\times\va}$ is transformed to $T_{0\times\va}$. Now from
the fact that $(n_i)$ has density~1 (see Remark~\ref{attractor})
and the fact that the action of $T_{0\times\va}$ preserves the set
$M_\mu(X\times G) \times\{\la_G\}$, it follows that each invariant
measure $P$ for $\wt{T}_{0\times\va}$ is concentrated on
$M_\mu(X\times G)\times\{\la_G\}$. Thus we only need to consider
the set of ergodic $(T\times Id_G)^{\wt{}}$-invariant measures on
$M_\mu(X\times G)$. Since the action of $T\times Id_G$  on
$X\times G$ is equicontinuous, we conclude that the system
$(C_2(\mu\ot\la_G),P,\wt{T}_{\va})$ has discrete spectrum.

Given $R\in C(T)$ denote by
$$
Z_R=\{\rho\in C_2(\mu\ot\la_G):\:\rho|_{X\times X}=\mu_R\}.$$
Notice that $Z_R$ is closed and $\wt{T_{\va}}$-invariant. By the
above we have proved the following result.

\begin{Prop}\label{FC102} If $T$ has discrete spectrum and $\va:X\to G$
is a cocycle taking values in a compact metric abelian group $G$ such that
$T_{\va}$ is ergodic then for each ergodic
$\wt{T}_{\va}$-invariant measure $P$ concentrated on $Z_{Id}$ the
resulting system $(Z_{Id},P,\wt{T}_{\va})$ has discrete
spectrum.
\bez\end{Prop}

\begin{Remark}\label{weier}\em
Let us note that the dificulty in the proof of Proposition~\ref{FC102}
arises from the fact that the cocycle $\va\times\va$ is not
ergodic. Indeed, we will shortly argue that whenever $T$ has
discrete spectrum and $\va:X \to G$ is ergodic
(i.e. $T_\va$ is ergodic) then the action of
$T_{\va}$ on $M_\mu(X\times G)$ has only one invariant measure
and moreover this measure is the Dirac measure supported on
$\mu\ot\la_G$.

First notice that by the Stone-Weierstrass theorem the set of
continuous functions of the form $f\ot\chi$, where $f\in C(X)$,
$\chi\in\widehat{G}$, is linearly
dense in $C(X\times G)$ and
therefore to check weak convergence of a sequence of measures in
$M_\mu(X\times G)$ we only need to test this on functions of the
form $f\ot\chi$. Assume that $(n_i)$ has density~$1$ and that
$U_{T_{\va}}^{n_i}$ converges (weakly) to zero in the orthocomplement
of $L^2(X,\mu)\ot1_G$. Take $\xi\in M_\mu(X\times G)$. Suppose
that along a subsequence of $(n_i)$, which we still denote by
$(n_i)$, we have $\left(T_{\va}^{n_i}\right)_\ast(\xi)\to\wt{\xi}$.
Then a calculation as above
shows that for $\chi\neq1$
$$
\int_{X\times G}f(x)\chi(g)\,d\wt{\xi}(x,g)=0,
$$
and we conclude that $\wt{\xi}=\mu\ot\la_G$.
\end{Remark}

Suppose now that the assumptions of Proposition~\ref{FC102} are
fulfilled and assume moreover that $R\in C(T)$ is such that the
cocycle
\beq\label{FC103}\va\times\va\circ R:X\to G\times
G\;\;\;\mbox{is ergodic.}
\eeq
We consider $\wt{T}_\va$ on $Z_R$.
Identyfying the two copies of $X$ we see that the action of
$\wt{T}_\va$ on $Z_R$ can be identified with a subsystem of the
action of $T_{\va\times\va\circ R}$ on $M_\mu(X\times G\times G)$.
In view of Remark~\ref{weier}, which by~(\ref{FC103}) we may apply
to $\va\times\va\circ R$, we obtain the following.

\begin{Prop}\label{FC104} Assume that $T$ has discrete spectrum
and that $\va:X\to G$ is a cocycle taking values in a compact metric
abelian group $G$ such that
$T_{\va}$ is ergodic. Then for each $R\in C(T)$ such that
$\va\times\va\circ R$ is ergodic the system $(Z_R,\wt{T}_\va)$ has
a quasi-attracting point. In particular it has only one invariant
measure (which is the Dirac measure concentrated at the unique
self-joining of $T_{\va}$ projecting on $\mu_R$).\bez\end{Prop}

Notice however that the two propositions above do not describe
all the ergodic invariant measures for $\wt{T}_\va$ in the case of $T$
having discrete spectrum. In particular we do not know whether an
ergodic measure $P$ (for $\wt{T}_\va$) which projects on the
natural Haar measure of $Y_J$, where the map $\kappa$
(see~(\ref{FC99})) is $1-1$, yields also a system with discrete
spectrum.

We will show next that going up one further step in the distal
echelon, namely taking a further isometric extension, one can already
obtain a transformation $\ov{T}$Ä preserving a measure $\ov{\mu}$
for which there exists an ergodic probability
measure $P \in M_{\wt{\ov{T}}}(C_2(\ov{\mu}))$ such that the
measure distal system
$(C_2(\ov{\mu}),P,{\wt{{\ov{T}}}})$ has mixed spectrum.

Let $Tx=x+\alpha$ be an irrational rotation on $\T$,
let $\va:\T\to\T$ be given by $\va(x)=x$ and let $T_{\va}:
\T^2 \to \T^2$ be defined by $T_{\va}(x,y)=(x+\alpha,x+y)$.
Let $Rx=x+\beta$, where
$\alpha,\beta$ and $1$ are independent over the rational numbers.
Then the cocycle
$\va\times\va\circ R$ is ergodic and the skew product
$T_{\va\times\va\circ R}$ has partly continuous spectrum (notice
however that its discrete part is bigger than the one given by the
rotation by $\alpha$). Now Proposition~\ref{FC104} applies and we
know that there exists only one (trivial) $\wt{T}_{\va}$-invariant
measure on $Z_R$.

Let now $\ov{T}$ (an extension
of $T_{\va}$) on $\T^3$ be given by
$$
\ov{T}(x,y,z)=(x+\alpha,x+y, x+y+z),
$$
and set $\ov{\mu} = \la_{\T}\ot\la_{\T}\ot\la_{\T}$.
Denote by $S_{a,b,c}\in
Aut(\ov{\mu})$ the rotation by
$(a,b,c)\in\T^3$. We have \beq\label{FC106} \ov{T}\circ
S_{a,b,c}\circ \ov{T}^{-1}=S_{a,a+b,a+b+c}.\eeq In other words the
action of $\wt{\ov{T}}$ on the set of graph measures given by
$S_{a,b,c}$, $(a,b,c)\in\T^3$, is isomorphic to the action $W$ on
$\T^3$ given by $$W(a,b,c)=(a,a+b,a+b+c).$$ Now, the homeomorphism
$W$ has many invariant tori $\T_a$ ($a\in\T$), where
$\T_a=\{(a,b,c):\:b,c\in\T\}$. Identifying $\T_a$ with $\T^2$, the
action of $W=W_a$ on $\T_a$ becomes $W_a(b,c)=(b+a,b+c+a)$. If now
$a=\beta$ is irrational then $W_\beta$ is isomorphic to the
$T_{\va}$ we started with (with $\beta$ taking the place
of $\al$; the constant cocycle $\beta$ is a
coboundary).  Thus we have proved the following.

\begin{Prop}\label{FC107}
For the transformation $\ov{T}(x,y,z)=(x+\alpha,x+y,x+y+z)$ as above
there are $\wt{\ov{T}}$-invariant measures $P$ such that the system
$(C_2(\la_{\T}^{\ot3}),P,\wt{\ov{T}})$ is ergodic and has partly
continuous spectrum. In fact, for each $\beta\in\T$ there exists
$P_\beta$ such that the system
$(C_2(\la_{\T}^{\ot3}),P_\beta,\wt{\ov{T}})$ is isomorphic to the
affine transformation $(x,y)\mapsto(x+\beta,x+y)$ of $\T^2$.
\bez\end{Prop}

We conclude with the following intriguing problem.

\begin{Problem}
Construct an ergodic system of zero entropy for which the
topological lens is a $P$-system. In view of Theorem
\ref{slabemiesz} this amounts to finding a weakly mixing system
$(X,\mathcal{B},\mu,T)$ with zero measure entropy such that the
periodic points of the topological lens are dense.
\end{Problem}


\begin{thebibliography}{99}




\bibitem{Co-Fo-Si}
I. P. Cornfeld, S. V. Fomin and Ya. G. Sinai,
{\em Ergodic theory\/},
Springer-Verlag, New York, 1982.

\bibitem{Da}
A. Danilenko,
{\em On simplicity concepts for ergodic actions},
J. Anal. Math. {\bfseries102}, (2007), 77-117.

\bibitem{EN}
R. Ellis and M. Nerurkar,
{\em Weakly almost periodic flows\/},
Trans.\ Amer.\ Math.\ Soc.\ {\bfseries 313}, (1989), 103-119.

\bibitem{Fa}
A. Fathi,
{\em Le groupe de transformations de $[0,1]$ qui preservent
la measure de Lebesgue est un groupe simple},
Israel J. of Math., {\bfseries 29}, (1978), 302-308.

\bibitem{Fu}
H. Furstenberg,
{\em Recurrence in ergodic theory and combinatorial number theory},
Princeton University Press,  Princeton, N.J., 1981.

\bibitem{Gl}
E. Glasner,
{\em Ergodic Theory via Joinings}, Mathematical Surveys and
Monographs {\bf 101}, AMS, Providence, RI, 2003.

\bibitem{Gl-Ki}
E. Glasner and J. King,
{\em A zero-one law for dynamical properties\/},
Topological dynamics and applications (Minneapolis,
MN, 1995), Contemporary Math.\  {\bfseries 215},
Amer. Math. Soc.,
(1998), 231-242.


\bibitem{Gl-We-93}
E. Glasner and B. Weiss,
{\em Sensitive dependence on initial conditions\/},
Nonlinearity  {\bfseries 6},
(1993), 1067-1075.


\bibitem{Ju-Ru}
A. del Junco and D.J. Rudolph,
{\em On ergodic actions whose
self-joinings are graphs\/}, Ergod.\ Th.\ Dynam.\ Sys.\ {\bfseries
7}, (1987), 531-557.

\bibitem{KR}
A. S. Kechris and C. Rosendal,
{\em Turbulence, amalgamation and
generic automorphisms of homogeneous structures\/},
Proc. Lond. Math. Soc. (3)  {\bfseries 94},  (2007),  302-350.

\bibitem{Kieff}
J. C. Kieffer,
{\em A simple development of the Thouvenot relative isomorphism theory},
Ann. Probab. {\bfseries12}, (1984),  204-211.


\bibitem{King}
J. King,
{\em  The generic transformation has roots of all orders},
Colloq. Math. {\bfseries 84/85}, (2000), part 2, 521-547.





\bibitem{Ku}
A. G. K\v{u}shnirenko, {\em Metric invariants of entropy type} (Russian),
Uspehi Mat. Nauk {\bfseries 22}, (1967) no. 5 (137), 57-65.

\bibitem{Le-Pa-Th}
M. Lema\'nczyk, F. Parreau and J.-P. Thouvenot,
{\em Gaussian automorphisms whose ergodic self-joinings
are Gaussian\/},
Fund.\ Math.\  {\bfseries 164},
(2000), 253-293.



\bibitem{OS}
D. Ornstein, {\em Bernoulli shifts with the same entropy are isomorphic},  Advances in Math.  {\bfseries 4},  (1970),  337-352.

\bibitem{O}
D. Ornstein,  {\em Ergodic theory, randomness, and dynamical systems},
James K. Whittemore Lectures in Mathematics given at Yale University. Yale Mathematical Monographs, No. 5. Yale University Press, 1974.


\bibitem{Th}
J.-P. Thouvenot,
{\em Quelques propri\'et\'es des syst\`emes dynamiques qui
se d\'ecomposent en un produit de deux syst\`mes dont l'un
est un sch\'ema de Bernoulli\/},
Israel J.\ of Math.\ {\bfseries 21}, (1975), 177-207.

\bibitem{Z1}
R. J. Zimmer,
{\em  Extensions of ergodic group actions\/},
Illinois J.\ Math.\ {\bfseries 20}, (1976), 373-409.

\bibitem{Z2}
R. J. Zimmer,
{\em  Ergodic actions with generalized discrete spectrum\/},
Illinois J.\ Math.\ {\bfseries 20}, (1976), 555-588.

\end{thebibliography}
\end{document}